\documentclass[number,sort&compress,preprint,11pt]{elsarticle}

\usepackage[a4paper,top=1in,bottom=1in,left=1.32in,right=1.32in]{geometry}

\usepackage{amsmath,amsthm,amsfonts,mathtools}

\usepackage{xspace}
\usepackage{enumitem}
\usepackage{breakurl}

\usepackage{tabularx}        
\usepackage{multirow}        
\usepackage{makecell}        
\usepackage{adjustbox}       
\usepackage[table,xcdraw]{xcolor} 
\usepackage{longtable}
\usepackage{array}
\usepackage{booktabs}

\newcolumntype{L}[1]{>{\raggedright\arraybackslash}m{#1}}
\newcolumntype{M}[1]{>{\centering\arraybackslash$\scriptsize}m{#1}<{$}}
\renewcommand{\arraystretch}{1.1}

\newcolumntype{L}[1]{>{\raggedright\arraybackslash}p{#1}}
\newcolumntype{C}[1]{>{\centering\arraybackslash}p{#1}}
\newcolumntype{R}[1]{>{\raggedleft\arraybackslash}p{#1}}

\usepackage{float}

\usepackage{graphicx}        
\usepackage{subcaption}
\usepackage{rotating}        
\usepackage{pdflscape}       

\usepackage{hyperref}

\newcommand{\setPhysicians}{\mathcal{P}}
\newcommand{\setPhysicianOneDutyWeekend}{\mathcal{P}^\textnormal{one}}
\newcommand{\setPhysicianSeveralDutiesWeekend}{\mathcal{P}^\textnormal{mult}}
\newcommand{\setManuallyPlannedPhysicians}{\mathcal{P}^\textnormal{manu}}
\newcommand{\setDuties}{\mathcal{D}}
\newcommand{\setMandatoryDuties}{\mathcal{D}^{\textnormal{mand}}}
\newcommand{\setManuallyPlannedDutyAssignments}{\mathcal{D}^{\textnormal{manu}}}
\newcommand{\setQualifiedDuties}{\mathcal{D}^{\textnormal{quali}}}
\newcommand{\setSoftQualifiedDuties}{\mathcal{D}^{\textnormal{quali}}_{\textnormal{soft}}}
\newcommand{\setQualifiedShift}{\mathcal{S}^{\textnormal{quali}}}
\newcommand{\setSoftQualifiedShift}{\mathcal{S}^{\textnormal{quali}}_{\textnormal{soft}}}
\newcommand{\setDays}{\mathcal{T}}
\newcommand{\setShifts}{\mathcal{S}}
\newcommand{\setDutyBlocks}{\mathcal{B}^\textnormal{duty}}
\newcommand{\setBlocksNoAddDuty}{\mathcal{B}^\textnormal{no-d}}
\newcommand{\setBlocksNoAddShift}{\mathcal{B}^\textnormal{no-s}}
\newcommand{\setShiftBlocks}{\mathcal{B}^\textnormal{shift}}
\newcommand{\setConsShiftBlocks}{\mathcal{B}^\textnormal{cons}}

\newcommand{\setMonths}{\mathcal{M}}
\newcommand{\setPools}{\Pi}
\newcommand{\setAbsences}{\mathcal{A}}
\newcommand{\setPhysicianPool}{\mathcal{P}}
\newcommand{\setPhysicianLastSB}{\mathcal{P}^\textnormal{prev}}
\newcommand{\setShiftPhysician}{\mathcal{P}}
\newcommand{\setManuallyPlannedShiftAssignments}{\mathcal{S}^\textnormal{manu}}
\newcommand{\setWeekend}{\mathcal{W}}
\newcommand{\setWeekendInMonth}{\mathcal{W}}
\newcommand{\setDutiesOfWeekend}{\mathcal{D}}
\newcommand{\setDutiesOfMonth}{\mathcal{D}}
\newcommand{\setDutiesOfPool}{\mathcal{D}}
\newcommand{\setDutyAbsences}{\mathcal{I}^\textnormal{duty}}
\newcommand{\setShiftAbsences}{\mathcal{I}^\textnormal{shift}}
\newcommand{\setDutyRestLastPP}{\mathcal{R}^\textnormal{duty}}
\newcommand{\setShiftRestLastPP}{\mathcal{R}^\textnormal{shift}}
\newcommand{\setDutySoftRestLastPP}{\mathcal{R}^\textnormal{duty}_\textnormal{soft}}
\newcommand{\setShiftSoftRestLastPP}{\mathcal{R}^\textnormal{shift}_\textnormal{soft}}
\newcommand{\setDutiesBeforeAbsence}{\mathcal{D}^\textnormal{before}}
\newcommand{\setDutiesAfterAbsence}{\mathcal{D}^\textnormal{after}}
\newcommand{\setConflicts}{\mathcal{C}}
\newcommand{\setConflictsSoft}{\mathcal{C}^\textnormal{soft}}

\newcommand{\parShiftMinimumNumber}{\textnormal{min}}
\newcommand{\parShiftDesiredNumber}{\textnormal{des}}
\newcommand{\parShiftMaxNumber}{\textnormal{max}}

\newcommand{\parRequiredFreeDaysAfterDutyBlock}{\textnormal{free\_days}}
\newcommand{\parStart}{\textnormal{start}}
\newcommand{\parEnd}{\textnormal{end}}
\newcommand{\parDesiredMaximumWeekends}{\textnormal{des\_max\_we}}
\newcommand{\parMaximumWeekends}{\textnormal{max\_we}}
\newcommand{\parPastConsecWeekends}{\textnormal{past\_we}}
\newcommand{\parDesiredMaximumConsecutiveWeekends}{\textnormal{cons\_we}}
\newcommand{\parMinFreeWeekends}{\textnormal{min\_free\_we}}
\newcommand{\parDesiredMinFreeWeekends}{\textnormal{des\_min\_free\_we}}
\newcommand{\parMaxNumberOfDutiesInPool}{\textnormal{max\_duties}}
\newcommand{\parDesiredMaxNumberOfDutiesInPool}{\textnormal{des\_max\_duties}}
\newcommand{\parMinNumberOfDutiesInPool}{\textnormal{min\_duties}}
\newcommand{\parDesiredMinNumberOfDutiesInPool}{\textnormal{des\_min\_duties}}
\newcommand{\parMaxPhysiciansPoolPerDay}{\textnormal{max\_phy}}
\newcommand{\parDesiredMaxPhysiciansPoolPerDay}{\textnormal{des\_max\_phy}}
\newcommand{\parExactDutiesPoolPerMonth}{\textnormal{ex\_pool}}
\newcommand{\parAttendencePhysician}{\textnormal{attend\_phy}}
\newcommand{\parFairDistribution}{\textnormal{target\_num}}

\newcommand{\parPrev}{\textnormal{prev}}
\newcommand{\parPhysicianLastD}{\textnormal{p\_prev\_pp}}
\newcommand{\parEndOfDays}{\textnormal{T}}

\newcommand{\MonthlyWeekendFactor}{\textnormal{we\_factor}}

\newcommand{\varDutyAssignment}{\textnormal{x}}
\newcommand{\varShiftAssignment}{\textnormal{y}}
\newcommand{\varDutyBlockAssignment}{\textnormal{x}^{\textnormal{block}}}
\newcommand{\varShiftBlockAssignment}{\textnormal{y}^{\textnormal{block}}}
\newcommand{\varShiftDesiredNumber}{\textnormal{y}^{\textnormal{desired}}}
\newcommand{\varShiftMaximalNumber}{\textnormal{y}^{\textnormal{max}}}
\newcommand{\varShiftDesiredNumberAuxiliary}{\textnormal{y}^{\textnormal{aux}}}
\newcommand{\varShiftBlockConsecPeriods}{\textnormal{y}^{\textnormal{block\_consec}}}
\newcommand{\varDutyConsecPeriods}{\textnormal{x}^{\textnormal{consec}}}
\newcommand{\varVioRestTime}{\textnormal{vio}^{\text{rest}}}
\newcommand{\varWeekendAttendence}{\textnormal{we}^{\text{attend}}}
\newcommand{\varWeekendViolation}{\textnormal{vio}^{\text{we\_pref}}}
\newcommand{\varMaxWeekendViolation}{\textnormal{vio}^{\text{max\_we}}}
\newcommand{\varMinFreeWeekendViolation}{\textnormal{vio}^{\text{free\_we}}}
\newcommand{\varMaxNumberOfDutiesViolation}{\textnormal{vio}^{\text{max\_d}}}
\newcommand{\varMinNumberOfDutiesViolation}{\textnormal{vio}^{\text{min\_d}}}
\newcommand{\varMaxNumberOfPhysiciansViolation}{\textnormal{vio}^{\text{max\_phy}}}
\newcommand{\varDeviationDutiesPoolNeg}{\textnormal{vio}^-}
\newcommand{\varDeviationDutiesPoolPos}{\textnormal{vio}^+}
\newcommand{\varMaxConsSBViolation}{\textnormal{vio}^{\text{max\_cons\_b}}}

\newcommand{\largeStepBackSmall}{\hspace{-0.75cm}}
\newcommand{\stepForLinkSizeSmall}{\hspace{0.6cm}}
\newcommand{\stepToGoBackToOriginSmall}{\hspace{0.1cm}}

\newcommand{\largeStepBack}{\hspace{-0.95cm}}
\newcommand{\stepForLinkSize}{\hspace{0.75cm}}
\newcommand{\stepToGoBackToOrigin}{\hspace{0.15cm}}

\begin{document}

\begin{frontmatter}
\title{A General Framework for Physician Rostering Using Mixed-Integer Programming and a Web-Based Graphical User Interface}

\author{Florian Meier\corref{mycorrespondingauthor}\fnref{tum}}
\ead{florian.meier@tum.de}
\author{Jan Boeckmann\fnref{linz}}
\ead{jan.boeckmann@jku.at}
\author{Clemens Thielen\fnref{tum}}
\ead{clemens.thielen@tum.de}
\cortext[mycorrespondingauthor]{Corresponding author}

\affiliation[tum]{organization={Technical University of Munich, Campus Straubing for Biotechnology and Sustainability, Professorship of Optimization and Sustainable Decision Making}, addressline={Am~Essigberg~3}, postcode={94315}, city={Straubing}, country={Germany}}
\affiliation[linz]{organization={Johannes Kepler University Linz, Institute for Business Analytics and Technology Transformation}, addressline={Altenberger~Straße~69}, postcode={4040}, city={Linz}, country={Austria}}

\begin{abstract}
Physician rostering in hospitals is complex due to varying shift structures, qualifications, and department- or hospital-specific regulations. Most existing optimization models are highly tailored to a single hospital or department and rarely see practical use. We present a general framework and a corresponding mixed-integer programming (MIP) model for physician rostering that accommodates a wide variety of roster structures and constraints. The model is integrated into a web application with an advanced graphical user interface (GUI), allowing physicians to specify preferences and hospital staff to configure the MIP model to their roster requirements without any mathematical or technical background. This approach enables easy adaptation to different hospitals or departments and straightforward updates in response to structural changes, such as new duties or modified qualifications. The applicability and effectiveness of the framework are demonstrated using real-world data from three departments in different hospitals specializing in internal medicine, cardiology, and orthopedics/trauma surgery. In one department, the system is already in everyday use, while in the other two, our model achieves comparable or improved roster quality relative to existing department-specific models, highlighting its potential as a versatile and practical tool for physician rostering.
\end{abstract}


\end{frontmatter}

\section{Introduction}\label{PR:sec:introduction}

Personnel planning is a central task in hospitals on all hierarchical levels --- from long-term strategic decisions concerning the required size and composition of the workforce to operational decisions about short-term adjustments of work schedules. \emph{Rostering problems} dealing with the design of work schedules (rosters) assigning duties or shifts to available personnel typically have planning horizons of a few weeks or months, and are therefore categorized as tactical or operational offline planning problems~\cite{erhard2018state}. These problems are mostly studied for nurses~\cite{cheang2003nurse,burke2004state,benazzouz2015literature} and physicians~\cite{erhard2018state} and have a large impact efficient hospital operations, quality of patient care, and employee satisfaction. Compared to nurse rostering, duty rostering problems for physicians are characterized by several specific features that make these problems particularly important, but also challenging to solve~\cite{erhard2018state}. While nurse rosters usually involve a fixed set of shift types with similar numbers of worked hours per day (e.g., an early shift, a late shift, and a night shift of eight hours each), physician rosters often contain more diverse sets of shift types that differ more substantially with respect to the number of hours worked in a shift or the required qualifications. Additionally, interactions of emergency duties at night and on weekends with assignments to specific workstations, clinical wards, or surgical teams during regular working hours often represent a crucial aspect to be considered in physician rostering problems~\cite{fugener2015duty,gross2018online,thielen2018duty}. Physicians are also a particularly valuable and expensive resource in hospitals and are considered a bottleneck in the care providing process~\cite{santos2014} and, in case they quit due to unsatisfactory working conditions, they are often hard to replace~\cite{Bodenheimer+Smith:primary-care}. This impacts the rostering process since it strongly increases the importance of aspects such as respecting the physicians' individual preferences regarding specific assignments and achieving a fair distribution of unpopular work at night or on weekends. Additionally, bad rosters that violate legally required limits on working times or collective labor contracts can lead to an increase in medical errors of physicians that can have severe consequences for patients.

\medskip

In practice, physician rosters in most hospitals in the world are still created manually by physicians who have to create complex rosters without any specific training or advanced digital decision support~\cite{Brunner+etal:flexible-shift,kraul2024optimizing}. Besides requiring a large amount of an expensive physician's valuable time, this often leads to unfair rosters that additionally violate labor regulations and do not adhere to the physicians' individual preferences~\cite{schoenfelder+pfefferlen:german-hospital}. Even though a large number of advanced models and algorithms for physician rostering exist in the Operations Research literature (see~\cite{erhard2018state} for a comprehensive overview), the number of cases where these methods are used in practice still remains extremely limited~\cite{brunner2011long,erhard2018state,thielen2018duty}. 

\medskip

An important reason for this is that the structure of physician rosters varies greatly across hospitals and departments, with differences in many key aspects~\cite{Rousseau+etal:general-approach}. Besides the number of physicians to consider or the length of the planning period, these aspects also include essential features defining how the rosters are organized. For instance, the number, working times, and required qualifications of duties and shifts usually differ significantly between institutions. Moreover, department- or physician-specific requirements such as individual labor contracts of physicians or templates that assign blocks of duties to the same physician for a certain time period such as a week often need to be respected. Since all of these requirements must be modeled in detail if a physician rostering method is to be used in practice, the majority of the existing approaches are designed for only a single hospital or department~\cite{erhard2018state}. The few cases where physician rostering models and algorithms from the literature have been successfully implemented in practice mostly involve a rostering method that has been specifically designed for an individual department of a given hospital that the researchers have collaborated with during the model's development. Due to the large variety of different requirements and roster structures, generic models that can be applied to a larger number of departments or hospitals would require, e.g., graphical user interfaces for entering individual preferences, labor agreements, and specific regulations of the hospital or department~\cite{Brunner+etal:flexible-shift,erhard2018state}.

\medskip

This paper contributes to addressing these issues by presenting a general problem framework and a corresponding mixed-integer programming (MIP) model for physician rostering that can accommodate a large variety of different requirements regarding the structure of the generated rosters. The model is implemented within a flexible web application with an advanced graphical user interface (GUI) that allows each physician to specify their individual preferences and view the generated rosters. In addition, the web application provides a comfortable way for hospital staff to configure the MIP model to a specific roster structure and specify all relevant parameters concerning the planned physicians and their specific qualifications as well as the duties and shifts to be assigned. Thereby, our rostering approach can be adapted to hospitals and departments with different roster structures without any changes to the MIP model or code. Moreover, in case of changes to the roster structure (e.g., due to newly introduced duties or changing qualification requirements of duties), our model can be easily adapted by hospital staff via the web application without any need for expert support.

We remark that our general framework and model can of course not represent every possible physician rostering problem that appears in the literature or in practice. Due to the extremely large variety in physician roster structures, this would be an extremely challenging or even impossible endeavor. Still, our model represents an important step towards more generic physician rostering approaches that can be flexibly adjusted to different hospitals or departments, which we demonstrate by testing it with real-word data from three departments that specialize in different disciplines (internal medicine, cardiology, and orthopedics/trauma surgery), are located in different hospitals, and use duty rosters with quite different structures. The considered cardiology department already uses our model and the web application in everyday practice. In the two other departments, department-specific MIP models for duty rostering have already been used for several years, and we compare the our general model to these department-specific models with respect to both the computation time and the quality of the generated rosters.

\subsection{Related Work}

Personnel planning problems for physicians in hospitals can be broadly categorized according to their planning horizons into \emph{staffing problems} (strategic decisions concerning the size and composition of the workforce), \emph{rostering problems} (tactical or operational offline creation of duty or shift rosters), and \emph{re-planning problems} (short-term, operational online adjustments of rosters or work schedules)~\cite{erhard2018state}. 
In the following, we focus on rostering problems. For an overview of the literature that includes staffing and re-planning problems for physicians, we refer to~\cite{erhard2018state}.

\medskip

Physician rostering is part of the larger field of personnel scheduling, which is a well-researched topic that has been the subject of numerous studies in the literature~\cite{van2013personnel,ozder2020systematic}. 
Regarding hospital settings, physician rostering is partly related to nurse rostering, which is reviewed, e.g., in~\cite{cheang2003nurse,burke2004state}. As outlined in the introduction and highlighted, e.g., in~\cite{erhard2018state,Brunner+etal:flexible-shift,thielen2018duty}, however, physician rostering exhibits a lot of properties that distinguish it from nurse rostering and lead to practical physician rostering problems being more complex and diverse in nature. In particular, the structure and constraints of physician rosters tend to differ much more across different departments and hospitals than is the case for nurse rosters~\cite{Brunner+etal:flexible-shift}.

\medskip

Numerous different aspects and extensions of physician rostering problems such as balancing the physicians' patient workloads~\cite{adams2019physician}, ensuring the viability of breaks~\cite{kraul2024optimizing}, or scheduling variable shift extensions to handle stochastic demand~\cite{fugener2019planning} have already been studied.
A comprehensive overview of the existing literature is provided in~\cite{erhard2018state}. As noted therein, however, only very few of the models presented in the literature have been successfully implemented in everyday practice, and most papers present approaches that are designed for only a single hospital or department.

\medskip

Physician rostering models and algorithms used in everyday practice include the emergency room physician scheduling methods presented in~\cite{carter2001scheduling}. They examine emergency room physician rosters from six hospitals and propose tabu search heuristics to optimize them, with two hospitals implementing the improved rosters in practice. Another example is the mixed-integer programming model described in~\cite{schoenfelder+pfefferlen:german-hospital}, which provides a detailed formulation of the rostering problem for an anesthesiology department in Germany. Their implementation includes an Excel-based interface that allows users to modify input data and adjust results in a straightforward way. A further application is presented in~\cite{fugener2015duty}, where two MIP models for assigning physicians to duties are developed and applied in a large German teaching hospital; a subsequent study~\cite{gross2018online} extends this work by proposing an online re-planning approach.

Another MIP approach applied in practice is presented in~\cite{thielen2018duty}. Here, a special focus is on the fair distribution of duties and the fulfillment of physicians' personal preferences, both of which play an important role in the overall satisfaction with the roster~\cite{bowers2016neonatal, dewa2017relationship}. A comparison of rosters generated by the MIP model from~\cite{thielen2018duty} to manually created ones over several planning periods demonstrates significant improvements not only in terms of fair distribution of duties and preference fulfillment, but also a significant reduction in understaffing of surgical teams during normal working hours caused by assignments of night and late duties. This clearly demonstrates the potential of using physician rostering models in practice and is in line with recent observations in~\cite{abdullah2025advancing}, where physician rosters generated by an MIP model are shown to yield an improvement of over~60\% in the objective value compared to manually-created rosters. Both papers~\cite{thielen2018duty,abdullah2025advancing}, however, state that their approaches are highly tailored to their specific partner hospitals and departments, which indeed applies to the far majority of models and algorithms presented in the literature.

\medskip

A notable exception is the already mentioned tabu search heuristics presented in~\cite{carter2001scheduling}, where emergency room physician rosters from six hospitals are examined. Furthermore, a case where a physician rostering model developed for a specific hospital is applied to additional instances that do not stem from this hospital is~\cite{cappanera2021emergency}. Here, the authors consider a two-phase planning problem where weekend shifts are first assigned to physicians over a six-month planning horizon in the first phase. The second phase then determines the workday shifts of each physician month by month while taking the weekend shifts assigned in the first phase into account. Besides real-world instances from a collaborating hospital, they also apply their MIP models to a set of instances derived from~\cite{wickert2021integer}.

The only previous papers we are aware of that explicitly aim to develop general models for physician rostering are~\cite{Rousseau+etal:general-approach,Bruni+Detti:physician-scheduling}. In~\cite{Rousseau+etal:general-approach}, a constraint programming (CP) model is proposed that introduces generic constraints for modeling the distribution of shifts 
and defining forbidden shift patterns. 
The model is implemented in ILOG OPL Studio Pro and tested on data from two hospital departments. While the generic definitions of distribution and pattern constraints allow the representation of various roster structures, the case studies involve only six classes of constraints, and these must be explicitly reformulated for each application using ILOG OPL Studio’s high-level modeling language. Hence, the approach cannot be directly applied to other departments or hospitals without changes to the model formulation. In~\cite{Bruni+Detti:physician-scheduling}, a flexible baseline mixed-integer programming model is presented and tested on data from a single hospital department. Similarly to~\cite{Rousseau+etal:general-approach}, the model's flexibility largely stems from a relatively low level of detail, involving only eight classes of constraints. Applying it in settings with additional requirements or different roster structures would therefore likely require modifications to the MIP in the form of additional variables or constraints. Consequently, this model also cannot be directly applied to other departments or hospitals without extending the underlying MIP formulation.
\section{Problem Description}\label{sec:problem_description}

We first describe the structure of the general physician rostering problem we consider. Afterwards, Sections~\ref{subsec:KL}--\ref{subsec:HD} illustrate the specific rostering problems in the three departments (internal medicine, cardiology, and orthopedics and trauma surgery) of different hospitals that are used to demonstrate the capabilities of our model and test its performance, and how these fit into our general problem structure. 

\subsection{General Problem Structure}\label{subsec:general_problem}
This section provides a detailed description of the general physician rostering framework we consider. Since the framework is designed to accommodate different rostering problems across hospital departments, the way a particular aspect or requirement appears in a specific scenario often depends on parameter settings, and some aspects may be irrelevant in certain cases (see Sections~\ref{subsec:KL}--\ref{subsec:HD} for concrete examples). All parameters that tailor the general framework to a specific department’s rostering problem can be configured via our graphical user interface (GUI) as described in Section~\ref{sec:implementation}.

\medskip

The overall goal of the rostering process is to assign a set of physicians to duties and shifts over a given planning period. The roster must satisfy a range of constraints, including structural requirements (e.g., staffing levels and necessary qualifications) and legal regulations (e.g., mandatory rest times). Its quality is evaluated based on several criteria, such as compliance with soft constraints, fairness, and the extent to which physician preferences are fulfilled.

\paragraph{Planning period}
We are given a planning period for which a duty roster is to be designed. The length of the planning period is typically one month, but in principle, any length or number of days is possible. For the days within the planning period, we distinguish between work days (days from Monday to Friday) and days on weekends (Saturdays and Sundays). Additionally, any day may be designated as a public holiday, since these days are often subject to specific rostering rules.

\medskip

\paragraph{Physicians and their qualifications}
There is a set of physicians with attributes such as their individual employment rate (in percent). Moreover, for each physician, there is a known set of days within the planning period on which the physician is absent (due to vacation, conference attendance, etc.). In addition, some physicians are designated as being \emph{planned manually}. For these physicians, all duties and shifts are pre-assigned outside the rostering process and are treated as fixed assignments that cannot be modified. 
Each physician can have any subset of a set of possible qualifications. These qualifications can represent either specific skills (e.g., expertise in critical care medicine) or experience levels (e.g., a certain number of completed years of specialist training) and can be defined and assigned to physicians via the GUI.

\medskip

\paragraph{Duties and their working times}
The main decision to be made in the rostering problem concerns the assignments of \emph{duties} to physicians on each day of the planning period. Each duty denotes a special activity (e.g., a night duty or staffing the catheterization lab) that should be assigned to exactly one physician on each day within a specific set of weekdays on which it occurs. Moreover, it can be specified whether the duty occurs also on public holidays and/or days before public holidays or whether it only occurs on these days. Some duties might already be pre-assigned to certain physicians and cannot be reassigned. Of the remaining ones, some mandatory duties must be assigned on each day they occur, while others can optionally remain unassigned. In addition, a duty has corresponding working times (e.g., 07:00 to 15:30), which might only end of the next day in case of a night duty (in which case the working times would be, e.g., 16:00 to 07:30). In this case, we still mean the day the night duty starts when referring a day where the duty occurs. Usually, the working times of a duty are identical on all days on which it occurs. If the activity corresponding to a duty has different working times on certain weekdays, a separate duty must be defined for each distinct working-time configuration. For public holidays and the days immediately before or after them, however, different working times can be specified for the same duty without creating a separate duty.

\medskip

\paragraph{Qualification rules for duty assignments}
For each duty, specific qualification constraints can be defined. A \emph{required qualification} means that a physician \emph{must} have the qualification to be assigned to the duty, while an \emph{excluded qualification} means that a physician \emph{must not} have it. In contrast, \emph{desired} and \emph{undesired qualifications} indicate that a physician \emph{should} or \emph{should not} have the corresponding qualification if possible. Undesired qualifications can, for instance, be used to prevent assignments of physicians with higher qualification levels (e.g., board-certified specialists) to duties that are intended for physicians still in specialist training.

\medskip

\paragraph{Shifts}
An important aspect of the problem is the interaction between duty assignments and the work physicians perform during their regular working hours. To clarify this distinction, we define \emph{duties} and \emph{shifts}. As described above, a duty denotes a special activity (e.g., a night duty or coverage of the catheterization lab) that should be assigned to exactly one physician on each day it occurs. In contrast, a shift refers to working on a specific ward during regular hours and may require multiple physicians per day. For each shift, both a mandatory minimum and an allowed maximum number of assigned physicians can be specified. The mandatory minimum, which can never be violated, is typically low (e.g., at least one physician on a ward). In addition, a higher desired minimum staffing level can be defined, which acts as a soft constraint whose violation incurs a user-defined penalty. Similar to duties, shifts can also involve required, excluded, desired, and undesired qualifications, and some shifts might already be pre-assigned to certain physicians. While shifts typically take place during regular hours from Monday to Friday, their exact working times may still vary (e.g., one shift runs from 07:00 to 15:30, another from 07:30 to 16:30).

\paragraph{Interaction of duties and shifts via rest times}
The interaction of duties and shifts stems from required rest times. Here, each duty or shift has defined \emph{rest times} that a physician should have after being assigned to it, which may depend on the next duty or shift she\footnote{For ease of reading, we use female pronouns (she/her) throughout the paper to refer to physicians, but these are intended to include physicians of all genders.} is assigned to. Here, for each ordered pair of duties or shifts, both a mandatory rest time that must always be respected between the two, as well as different levels of desired rest times that can potentially be violated are possible. For instance, a mandatory rest time of 24~hours is sometimes required between two night duties (which would, e.g., prevent night duties by the same physician on Monday and Tuesday), but it might be desirable that a physician does not have several night duties in short succession (e.g., on Monday and Wednesday or Monday and Thursday). In this case, additional desired rest times of 48 hours (discouraging assignments on both Monday and Wednesday) and 72 hours (discouraging an assignment on both Monday and either Wednesday or Thursday) would be specified between any two night duties, where respecting the 48 hour rest time would be of higher importance. Note that, a positive rest time between two duties or shifts in particular means that these two must not be assigned to the same physician simultaneously (i.e., with overlapping working times). For cases where two assignments are allowed to partially overlap, a negative rest time can be used to specify the permitted number of hours of overlap.

\medskip

\paragraph{Effects of absences on duty and shift assignments}
In addition to rest times, the given absences of physicians must be respected. Clearly, no duty or shift can be assigned to a physician on a day when she is absent. While shifts may still be assigned on the days immediately before or after an absence, each duty specifies separately whether it may be assigned on the day preceding an absence and whether it may be assigned on the day following an absence. For instance, night duties are usually prohibited on the day before an absence.

\medskip

\paragraph{Duty blocks and shift blocks}
Another important aspect of the problem is that certain duties and shifts are grouped and assigned as \emph{blocks} rather than individually to ensure continuity in both physician working times and patient care. A \emph{duty block} (\emph{shift block}) is a set of duties (shifts) on specific days that must all be assigned to the same physician. Such blocks are typically defined over a period of about one week, although different time spans are also possible. For example, a duty block might consist of catheterization lab duties on Monday and Tuesday and a night duty on Friday and Saturday, all of which must be assigned to the same physician. For some duty or shift blocks, the assigned physician may also take on additional duties or shifts in between, provided that all required rest times are respected, whereas for other blocks this is not permitted. Moreover, a duty or shift block can require a block-specific number of free days after the block. Note that either the block itself or the free days required after it may extend into the following planning period, and these constraints must be taken into account when creating the next roster.

\medskip

\paragraph{Consecutive assignments}
In addition to predefined duty and shift blocks, continuity of working times and patient care can also be promoted by specifying desired consecutive assignments of duties or shift blocks. For a duty, this means that it is desirable (but not mandatory) to assign it to the same physician on several consecutive days --- for instance, to avoid frequent switches between a late duty and regular early-morning shifts. Similarly, it may be desirable to assign certain shift blocks to the same physician(s) across multiple consecutive assignments to promote continuity of care on the corresponding ward. At the same time, an upper bound can be imposed on the number of consecutive assignments of the same shift block to prevent a physician from staying on the same ward for too long (e.g., to ensure sufficient variety in their tasks for training purposes).

\medskip

\paragraph{Fair distribution of duties}
A major goal of the rostering process is to distribute duties fairly among physicians. There exists a wide variety of fairness concepts, metrics, and implementations in the physician rostering literature (see~\cite{fuchs2025fairness} for a recent overview). Approaches to fairness can also differ significantly between hospitals and departments. To accommodate these differences, we consider a general fairness model that can represent various implementations of fairness.
To this end, we introduce the concept of a \emph{pool}, defined by a set of duties (the \emph{duties in the pool}) and a set of physicians (the \emph{physicians in the pool}). Pools can be created in the GUI, and each pool allows the specification of hard and/or soft lower and upper bounds on the number of duties in the pool each physician in the pool may be assigned. Alternatively, an exact number of duties in the pool that each physician in the pool must be assigned can be specified. 
Additionally, a fair distribution of the duties among the physicians in the pool can be specified. In this case, it is desirable to assign each physician a number of duties proportional to their availability, measured as the number of days they are not absent multiplied by their employment rate. The employment rate is only factored in when part-time employment does not simply result in additional absent days, since these are already accounted for. Taking absences into account ensures that physicians who are away for extended durations (e.g., on holiday) do not receive an unusually large number of duties for the remainder of the planning period.
Also note that the sets of physicians and duties in different pools do not need to be disjoint. For example, one pool may include a subset of physicians and all night duties, ensuring lower and upper bounds and a fair distribution of all night duties for these physicians, while another pool with the same physicians but only the night duties on Saturdays enforces similar constraints specifically for Saturday night duties. In this way, very different implementations of fairness can be represented using pools.

\medskip

\paragraph{Weekends}
An important aspect that also relates to fairness is the distribution of work on weekends, which is often regulated by collective labor agreements. Here, a weekend is considered a \emph{worked weekend} for a physician if she is assigned to any duty on Saturday or Sunday, or to a duty on Friday whose working time ends after a specified threshold (e.g., 21:00). Otherwise, the weekend is considered a \emph{free weekend}. Possible requirements regarding weekends include hard upper bounds on the number of consecutive worked weekends, as well as hard and/or soft upper bounds on the number of worked weekends per month. Alternatively, hard and/or soft lower bounds on the number of free weekends per month are possible as some collective labor agreements require such lower bounds explicitly. Note that, in addition to these bounds on worked or free weekends, pools can be used to foster a fair distribution of duties on weekends overall or of specific duties on weekends that are particularly undesired by most physicians (e.g., night duties on Saturdays). These pools relating to weekend duties often include also duties on or before public holidays.

\medskip

\paragraph{Simultaneous assignments}
Besides modeling fairness, pools are also used to represent operational constraints concerning simultaneous assignments of physicians from specific groups to duties (e.g., that at most two physicians with a low experience level can be assigned to concurrent night duties). To this end, hard and/or soft upper bounds on the number of physicians in a pool that can be assigned to duties in this pool on the same day can be specified. 

\paragraph{Physician preferences}
Besides achieving a fair roster, respecting the physicians' individual preferences is also a crucial goal of the rostering process. To this end, each physician can specify her preferences for the planning period via the GUI. Like the specific approaches for handling fairness, the ways in which physicians can express their preferences often differ significantly between departments and hospitals. Our general preference model distinguishes three main types of physician preferences, not all of which need to be applied in every scenario: (i) \emph{duty-specific preferences}, which refer to a particular duty on a given day; (ii) \emph{weekly preferences}, which apply to a user-defined set of duties and/or shifts within a week; and (iii) \emph{weekend preferences}, which concern the assignment of duties on worked weekends.
Providing duty-specific preferences might only be allowed on certain days of the week. The selectable options can include \emph{strongly desired}, \emph{desired}, \emph{indifferent}, or \emph{undesired}, which do not necessarily have to be satisfied in the generated roster but can be weighted by their relative importance. Additionally, an option \emph{impossible}, which must be strictly enforced, can also be offered. For each preference option, an upper bound on the number of times a physician may select it can be specified.
Weekly preferences allow arbitrary combinations of duties and shifts within a week (e.g., the duty consisting of staffing the catheterization lab from Monday to Thursday and all regular shifts on a specific ward from Monday to Friday) to be grouped into a single preference set. Each physician qualified for any of the corresponding duties and shifts can then specify, for each week of the planning period, how desirable the entire set is for her. The available options are the same as for duty-specific preferences, except that \emph{impossible} cannot be used for weekly preferences.
Weekend preferences are specified once per physician for the entire planning period. They indicate whether the physician prefers fewer worked weekends with multiple duties assigned on each, or a larger number of worked weekends with only one duty assigned per weekend.

\medskip

\paragraph{Effects of assignments from the previous planning period}
While duties and shifts are only assigned for the days of the current planing period, some of the described constraints on the duty roster require that assignments made in the previous planning period are taken into account. This is relevant to ensure rest times or required free days between previously assigned duties, shifts, and blocks and assignments made in the current planning period. Moreover, the previous assignments also need to be considered when evaluating consecutive assignments of duties, shift blocks, and worked weekends. Furthermore, when the current planing period does not start at the beginning of a month, previous assignments are also relevant when evaluating the number of worked weekends within a month.

\subsection{Problem Description for Internal Medicine Department}\label{subsec:KL}
The internal medicine department we consider has 128 beds, including a 15-bed intensive care unit, and around 35 physicians, some working part-time (60\% or 80\%). The planning period is typically one calendar month, sometimes including adjacent public holidays (e.g., New Year’s Day in the December planning period).

\paragraph{Duties and shifts}
Duties comprise two night duties (N1, N2; 20:00--08:00), which occur every day, and two day duties (D1, D2; 08:00--20:00), which occur only on weekends and public holidays. All regular duties are mandatory. Each duty also has a corresponding \emph{backup duty}, which is optional to staff: the physician assigned to a backup duty does not normally work during the scheduled hours, but steps in if the originally assigned physician becomes unavailable on short notice. Certain duties require specific qualifications (e.g., six months of intensive care experience), while some physicians are excluded from particular or all duties due to contracts or health reasons. Physicians are assigned to fixed wards for regular shifts (07:15--16:00, Monday to Friday), which is represented by required qualifications. Minimum staffing levels for ward shifts ensure adequate coverage, with no upper bound.

\medskip

\paragraph{Interaction via rest times}
At most one regular duty or shift can be assigned to a physician per day, and mandatory rest times between duties and shifts ensure sufficient recovery. Desired rest times are used to discourage night duties with only one or two duty-free days in between. Backup duties follow the same rules, except that a backup duty may be assigned on the same day as a shift.

\medskip

\paragraph{Absences}
Planned absences prevent duty, backup duty, or shift assignments on absent days. Duties are also prohibited on the day before an absence, while shifts remain allowed.

\medskip

\paragraph{Duty and shift blocks}
In this department, all regular and backup duties are assigned individually, so no duty blocks are used. Likewise, since each physician is assigned to a fixed ward, shift blocks are not required. Consequently, consecutive assignments of shift blocks do not arise, and consecutive duties are prevented by the mandatory rest times described above.

\medskip

\paragraph{Fair distribution via pools}
Some physicians have fixed reduced duty counts due to contracts or health reasons. Remaining duties are distributed fairly using \emph{pools}: each physician’s monthly target number of duties (excluding backups) is proportional to their availability (non-absent days times employment rate). Backup duties are excluded from fairness constraints but limited to four per month. Additional pools enforce weekend rules (e.g., limiting Saturday night duties to once every two months) and simultaneous assignment constraints (e.g., at most one physician from each ward section on night duty per day).

\medskip

\paragraph{Physician preferences}
Only duty-specific preferences are used, with the options \emph{strongly desired}, \emph{desired}, \emph{indifferent}, \emph{undesired}, and \emph{impossible}. The option \emph{undesired} does not strictly prevent assignments and may be selected up to ten times per month per physician, whereas \emph{impossible} strictly prohibits any assignment on that day and may be used on up to three days per month.

\medskip

\paragraph{Assignments from the previous period}
Assignments from the previous planning period are considered to enforce rest times and limit consecutive night duties at the start of the current period.

\subsection{Problem Description for Cardiology Department}\label{subsec:SR}
The cardiology department has about 30 physicians, who staff several regular wards, an intensive care unit (ICU), and a chest pain unit (CPU). The planning period is typically one month, sometimes adjusted due to public holidays.

\paragraph{Duties and shifts}
Duties include intermediate duties (15:30--00:00), night duties (22:00--08:00), ICU duties (early, day, late, night), CPU duties (early, late, night, outpatient), a Function duty, which is responsible for performing electrocardiography and echocardiography during regular working hours from Monday to Thursday, and a Function Support duty, which covers the functional area of wards from other departments from Monday to Friday. On weekdays, most duties are organized in weekly \emph{duty blocks} and shifts in \emph{shift blocks}. Some duty blocks even combine multiple types of duties, e.g., a block may include the Function duty from Monday to Thursday and a night duty on Friday and Saturday. On weekends and public holidays, duties are assigned individually. Some duties require specific qualifications (ICU or CPU), while others have no qualification requirement. For certain shifts on specific wards, it is also desired that only physicians with the corresponding qualifications are assigned. 
Because the department faces severe physician shortages, all duties are modeled in our framework as \emph{optionally unassigned}. It is often impossible to staff every duty with the department’s own personnel, so unfilled duties are typically covered by physicians from other departments or by costly external staff. Similarly, most minimum staffing levels for ward shifts are represented in our framework as \emph{desired minimum staffing levels} (soft constraints) that can be violated if necessary. Consequently, assigning as many duties as possible and respecting these soft constraints is a central planning objective.

\medskip

\paragraph{Interaction via rest times}
Each physician can be assigned to at most one duty or shift per day. Mandatory rest times ensure sufficient recovery between consecutive assignments, while pairs of duties that should not be assigned to the same physician if possible are represented as desired rest times in our framework.

\medskip

\paragraph{Absences}
Physicians’ planned absences are known prior to the planning process. A physician cannot be assigned to any duty or shift on a day they are absent. In addition, duties that end after 21:00 cannot be assigned to a physician with an absence on the following day.

\medskip

\paragraph{Fair distribution via pools}
To ensure fairness, most duties are assigned using \emph{pools} in our framework. This includes fair distribution of ICU and CPU duties, weekend duties, and other demanding or undesired duties, taking into account individual pre-assigned duties and contractual exceptions. Physicians may also work at most five Function Support duties per month. Upper bounds on pools are additionally used to model operational constraints for \emph{simultaneous assignments}, e.g., limiting the number of physicians from a specific group assigned to certain duties on the same day, such as the restriction that at most one of ICU early and ICU day duty on the same day may be assigned to a physician marked as new on the ICU.

\medskip

\paragraph{Physician preferences}
Duty-specific preferences are allowed on weekends and public holidays, with the options \emph{strongly desired}, \emph{desired}, \emph{undesired}, \emph{indifferent}, and \emph{impossible} (at most two times per planning period). During the week, weekly preferences apply instead: duties and shifts are categorized into morning, day, evening, and night, and physicians can indicate a preference (\emph{strongly desired}, \emph{desired}, \emph{undesired}, or \emph{indifferent}) for each category. Weekend preferences indicate whether a physician prefers fewer worked weekends with multiple duties or more worked weekends with only one duty.

\medskip

\paragraph{Assignments from the previous period}
Assignments from the previous planning period are considered to enforce rest times, prevent consecutive night duties at the start of the current period, and account for duty and shift blocks that span the boundary between the previous and current planning periods. Rest times and continuity constraints arising from blocks or duties that started in the previous period are therefore taken into account when assigning duties and shifts in the current period.

\subsection{Problem Description for Orthopedics/Trauma Surgery Department}\label{subsec:HD}
The orthopedics/trauma surgery department’s duty rostering problem has been described in detail in~\cite{thielen2018duty}. Since that publication, some structural changes have been made (e.g., the removal of float physicians and the introduction of a new collective agreement governing fairness). The department has about 50~physicians, including residents and fellows. The planning period is typically one month, with slight adjustments in some cases due to public holidays.

\paragraph{Duties and shifts}
Duties include five night duties and one late duty per day. Night duties start after regular working hours on weekdays and in the morning on weekends and public holidays, ending the following morning. The late duty runs from mid-afternoon to late evening, with slightly earlier start and end times on weekends and public holidays. All duties are mandatory to staff and at most one of them can be assigned to each physician per day. Moreover, each duty requires a duty-specific experience level, and some physicians are excluded from particular duties due to contractual agreements.
During regular working hours, each physician is permanently assigned to one of 12 surgical wards for regular shifts from 07:15 to 16:00 on weekdays (excluding public holidays). These permanent assignments are represented by required qualifications for the shifts of each ward.

\paragraph{Consecutive assignments}
For the late duty, it is desirable, if possible, to assign the same physician to multiple consecutive occurrences to promote continuity of patient care and working times. This is modeled in our framework as a desired consecutive assignment.

\paragraph{Interaction via rest times}
A key challenge in this department is the interaction between duty assignments and surgical ward staffing. Each ward has daily minimum staffing requirements for both the total number of physicians and the number of fellows. In our framework, these are represented as desired minimum staffing levels, i.e., soft constraints that may be violated in exceptional cases (e.g., if too many physicians of a ward are absent due to vacation). Duty assignments affect these constraints through their interaction via mandatory rest times: physicians are unavailable on their ward on the day of a late duty and on the day after a night duty (except for a stand-by night duty), but remain available on the day of a night duty since they first work their ward shift before starting the duty. In addition, pre-assigned ward duties and planned absences are modeled as restrictions preventing duty assignments on the affected days.

\medskip

\paragraph{Fair distribution via pools}
Following a new collective agreement, which has replaced the detailed fairness model described in~\cite{thielen2018duty}, fairness is now mainly enforced through upper bounds: each physician may work at most four night duties per month and at most two weekends per month (where a weekend counts if any duty is assigned on Friday, Saturday, or Sunday). Pools in our framework represent these constraints and ensure that violations, which trigger monetary compensation, are used only in exceptional cases.

\medskip

\paragraph{Physician preferences}
Duty-specific preferences are allowed for all duties on days when physicians are neither absent nor pre-assigned to their ward. The options are \emph{desired}, \emph{indifferent}, and \emph{undesired}. The option \emph{undesired} can be used on at most half of the days in the planning period, and on at most half of the weekend days, including public holidays.

\medskip

\paragraph{Assignments from the previous period}
Assignments from the previous planning period are considered to enforce rest times and prevent consecutive night duties at the start of the current period.

\section{Mathematical Modelling}\label{PR:sec:MIP}

We now present the detailed mixed-integer programming (MIP) formulation that models the problem described in Section~\ref{subsec:general_problem}. We begin by listing all sets, parameters, and variables in Tables~\ref{tab:sets}--\ref{tab:variables}, respectively. Subsequently, we explain the individual components of the objective function and describe the model constraints. While the constraints are introduced verbally in the main text, their exact mathematical formulations are provided in \ref{PR:sec:appendix_MIP_formulation}.

\medskip

As will be discussed in Section~\ref{sec:implementation}, the model can be configured to represent a specific rostering problem---such as the three examples described in Sections~\ref{subsec:KL}--\ref{subsec:HD}---via our graphical user interface. Consequently, the model generated for a given rostering problem includes only those parts that are relevant to the corresponding application. For example, variables and constraints related to duty and shift blocks are omitted in cases such as the internal medicine department, where no such blocks are used (see Section~\ref{subsec:KL}). This modular structure ensures that, despite the high flexibility of our general formulation, each instantiated MIP remains as compact and efficient as possible.

\begin{longtable}{|p{0.14\linewidth}|p{0.79\linewidth}|}
\caption{Model sets with descriptions.\label{tab:sets}} \\
\hline
\textbf{Set} & \textbf{Description} \\
\hline
\endfirsthead

\multicolumn{2}{c}%
{\tablename\ \thetable\ -- \textit{continued from previous page}}\\
\hline
\textbf{Set} & \textbf{Description} \\
\hline
\endhead

\hline \multicolumn{2}{r}{\textit{continued on next page}} \\
\endfoot
\endlastfoot %

$\setDays = \{1, \dots , \parEndOfDays\}$ & Set of days (index~$t$) \\ \hline
$\setWeekend = \{1, \dots , W\}$ & Set of weekends (index~$w$) \\ \hline
$\setMonths$ & Set of months (index~$m$) \\ \hline
$\setWeekendInMonth(m)$ & Set of weekends in month~$m$ \\ \hline
$\setPools$ & Set of pools (index~$\pi$) \\ \hline
$\setPhysicians$ & Set of physicians (index~$p$) \\ \hline
$\setManuallyPlannedPhysicians$ & Set of physicians that are planned manually \\ \hline
$\setPhysicianPool(\pi)$ & Set of physicians in pool~$\pi$ \\ \hline
$\setShiftPhysician(s)$ & Set of physicians in ward of shift~$s$ \\ \hline
$\setPhysicianOneDutyWeekend$ & Set of physicians who prefer to have one duty per weekend but more worked weekends \\ \hline
$\setPhysicianSeveralDutiesWeekend$ & Set of physicians who prefer to have multiple duties per weekend but fewer worked weekends \\ \hline
$\setPhysicianLastSB(b)$ & Set of physicians who should be assigned to shift block~$b$ due to their assignment to the last shift block of the previous planning period, intended for consecutive assignment \\ \hline
$\setDuties$ & Set of duties (index~$d$) \\ \hline
$\setMandatoryDuties$ & Set of mandatory duties \\ \hline
$\setManuallyPlannedDutyAssignments(p)$ & Set of duties to which physician~$p$ is manually assigned \\ \hline
$\setQualifiedDuties(p)$ & Set of duties for which physician~$p$ has all required qualifications and does not have any of the excluded qualifications \\ \hline
$\setSoftQualifiedDuties(p)$ & Set of duties for which physician~$p$ has all required and desired qualifications and does not have any of the undesired or excluded qualifications \\ \hline
$\setDuties(t)$ & Set of duties on day~$t$ \\ \hline
$\setDutiesOfWeekend(w)$ & Set of duties on weekend~$w$ \\ \hline
$\setDutiesOfMonth(m)$ & Set of duties in month $m$ \\ \hline
$\setDutiesOfPool(\pi)$ & Set of duties in pool $\pi$ \\ \hline
$\setDutiesBeforeAbsence$ & Set of duties that cannot be assigned to a physician on a day before an absence \\ \hline
$\setDutiesAfterAbsence$ & Set of duties that cannot be assigned to a physician on a day after an absence \\ \hline
$\setShifts$ & Set of shifts (index~$s$) \\ \hline
$\setQualifiedShift(p)$ & Set of shifts for which physician~$p$ has all required qualifications and does not have any of the excluded qualifications \\ \hline
$\setSoftQualifiedShift(p)$ & Set of shifts for which physician~$p$ has all required and desired qualifications and does not have any of the undesired or excluded qualifications \\ \hline
$\setShifts(t)$ & Set of shifts on day~$t$ \\ \hline
$\setManuallyPlannedShiftAssignments(p)$ & Set of shifts to which physicians~$p$ is manually assigned \\ \hline
$\setDutyBlocks$ & Set of duty blocks (index~$b$) \\ \hline
$\setShiftBlocks$ & Set of shift blocks (index~$b$) \\ \hline
$\setBlocksNoAddDuty$ & Set of duty and shift blocks during which a physician cannot be assigned additional duties \\ \hline
$\setBlocksNoAddShift$ & Set of duty and shift blocks during which a physician cannot be assigned additional shifts \\ \hline
$\setConsShiftBlocks = \{\setConsShiftBlocks_1,\allowbreak \setConsShiftBlocks_2,\allowbreak \dots\}$
& Set containing all sets $\setConsShiftBlocks_{j} = \{b^j_{1}, b^j_{2}, \dots\}$ of consecutive shift blocks that should not all be assigned to the same physician due to a desired maximum number of consecutive shift blocks. If no physician should be assigned more than $n$~consecutive shift blocks, then the corresponding set $\setConsShiftBlocks_{j}$ contains $n+1$~consecutive shift blocks \\ \hline
$\setAbsences(p)$ & Set of days on which physician~$p$ is absent \\ \hline
$\setDutyAbsences(p)$ & Set of duties that cannot be assigned to physician~$p$ since they have selected \emph{impossible} as their preference  \\ \hline
$\setShiftAbsences(p)$ & Set of shifts  that cannot be assigned to physician~$p$ since they have selected \emph{impossible} as their preference  \\ \hline
$\setConflicts$ & Set of pairs~$(d_1, d_2)$ of duties and shifts ($d_1, d_2 \in\setDuties\cup\setShifts$) that cannot be assigned to the same physician due to mandatory rest times \\ \hline
$\setConflictsSoft$ & Set of pairs~$(d_1, d_2)$ of duties and shifts ($d_1, d_2 \in\setDuties\cup\setShifts$) that should not be assigned to the same physician due to desired rest times \\ \hline
$\setDutyRestLastPP(p)$ & Set of duties to which physician~$p$ cannot be assigned due to mandatory rest times after assigned duties or shifts of the previous planning period \\ \hline
$\setShiftRestLastPP(p)$ & Set of shifts to which physician~$p$ cannot be assigned due to mandatory rest times after assigned duties or shifts of the previous planning period \\ \hline
$\setDutySoftRestLastPP(p)$ & Set of duties to which physician~$p$ should not be assigned due to desired rest times after assigned duties or shifts of the previous planning period \\ \hline
$\setShiftSoftRestLastPP(p)$ & Set of shifts to which physician~$p$ should not be assigned due to desired rest times after assigned duties or shifts of the previous planning period \\ \hline
\end{longtable}

\begin{longtable}{|p{0.25\linewidth}|p{0.71\linewidth}|}
\caption{Model parameters with descriptions.\label{tab:parameters}} \\

\hline
\textbf{Parameter} & \textbf{Description} \\
\hline
\endfirsthead

\multicolumn{2}{c}%
{\tablename\ \thetable\ -- \textit{continued from previous page}}\\
\hline
\textbf{Parameter} & \textbf{Description} \\
\hline
\endhead

\hline \multicolumn{2}{r}{\textit{continued on next page}} \\
\endfoot

\endlastfoot

$\parShiftMinimumNumber(s)$ & required minimum number of physicians assigned to shift~$s$ \\ \hline
$\parShiftDesiredNumber(s)$ & desired minimum number of physicians assigned to shift~$s$ \\ \hline
$\parShiftMaxNumber(s)$ & maximum allowed number of physicians assigned to shift~$s$ \\ \hline
$\parRequiredFreeDaysAfterDutyBlock(b)$ & required number of free days after the duty or shift block~$b$ \\ \hline
$\MonthlyWeekendFactor(m)$ & fraction of Saturdays of month~$m$ that are in the current planning period \\ \hline
$\parMaximumWeekends$ & maximum allowed number of worked weekends per physician and month \\ \hline
$\parDesiredMaximumWeekends$ & desired maximum number of worked weekends per physician and month \\ \hline
$\parPastConsecWeekends(p)$ & number of consecutive weekends physician~$p$ worked at the end of the previous planning period \\ \hline
$\parDesiredMaximumConsecutiveWeekends$ & maximum number of consecutive weekends a physician is allowed to work \\ \hline
$\parMinFreeWeekends$ & minimum required number of free weekends per physician and month \\ \hline
$\parDesiredMinFreeWeekends$ & desired minimum number of free weekends per physician and month \\ \hline
$\parMaxNumberOfDutiesInPool(\pi)$ & maximum allowed number of duties in pool~$\pi$ for each physician in the pool. This parameter is set to \texttt{None} if no maximum number is specified \\ \hline
$\parDesiredMaxNumberOfDutiesInPool(\pi)$ & desired maximum number of duties in pool~$\pi$ for each physician in the pool. This parameter is set to \texttt{None} if no desired maximum number is specified \\ \hline
$\parMinNumberOfDutiesInPool(\pi)$ & required minimum number of duties in pool~$\pi$ for each physician in the pool. This parameter is set to \texttt{None} if no minimum number is specified \\ \hline
$\parDesiredMinNumberOfDutiesInPool(\pi)$ & desired minimum number of duties in pool~$\pi$ for each physician in the pool. This parameter is set to \texttt{None} if no desired minimum number is specified \\ \hline
$\parMaxPhysiciansPoolPerDay(\pi)$ & maximum allowed number of physicians in pool~$\pi$ assigned to duties in the pool on the same day. This parameter is set to \texttt{None} if no maximum number is specified \\ \hline
$\parDesiredMaxPhysiciansPoolPerDay(\pi)$ & desired maximum number of physicians in pool~$\pi$ assigned to duties in the pool on the same day. This parameter is set to \texttt{None} if no desired maximum number is specified \\ \hline
$\parExactDutiesPoolPerMonth(\pi)$ & exact number of duties in pool~$\pi$ that must be assigned to each physician in the pool. This parameter is set to \texttt{None} if no exact number is specified \\ \hline
$\parFairDistribution(\pi, p)$ & fair number of duties in pool~$\pi$ that should be assigned to physician~$p$ according to her volume of employment, number of absent days, and other physicians in the pool. This parameter is set to \texttt{None} if no fair distribution is required in the pool \\ \hline
$\parPrev(b)$ & preceding shift block of shift block~$b$ that should be assigned to the same physicians. This parameter is set to \texttt{None} if no previous shift block should be assigned to the same physicians \\ \hline
$\parPrev(d)$ & preceding duty of duty~$d$ that should be assigned to the same physician. This parameter is set to \texttt{None} if no previous duty should be assigned to the same physician \\ \hline
$\parStart(b)$ & first day of shift block or duty block~$b$ \\ \hline
$\parEnd(b)$ & last day of shift block or duty block~$b$ \\ \hline
$\parPhysicianLastD(d)$ & physician assigned to the last duty of the previous planning period who should also be assigned to duty~$d$. This parameter is set to \texttt{None} if no physician is assigned to $\parPrev(d)$ or $\parPrev(d)$ is \texttt{None} \\ \hline
\end{longtable}

\medskip
For each pool~$\pi\in\setPools$ and each physician~$p\in\setPhysicians(\pi)$, the target number $\parFairDistribution(\pi, p)$ is calculated as
\begin{align*}
    \parFairDistribution(\pi, p) &= n(\pi) \cdot\frac{ v(p) \cdot \parAttendencePhysician(p)}{\sum_{p\in\setPhysicianPool(\pi)} v(p) \cdot \parAttendencePhysician(p)},
\end{align*}
where $n(\pi)$ denotes the number of duties in pool~$\pi$, $v(p)$ the employment volume of physician~$p$, and $\parAttendencePhysician(p)$ the number of duties on whose corresponding days physician~$p$ is not absent.

\begin{longtable}{|p{0.16\linewidth}|p{0.61\linewidth}|p{0.14\linewidth}|}
\caption{Model variables with descriptions and objective function coefficients. Positive coefficients encourage (i.e., reward) assignments, while negative coefficients penalize them.\label{tab:variables}}\\
\hline
\textbf{Variable} & \textbf{Description} & \textbf{Objective weight} \\
\hline
\endfirsthead

\multicolumn{3}{c}%
{\tablename\ \thetable\ -- \textit{continued from previous page}}\\
\hline
\textbf{Variable} & \textbf{Description} & \textbf{Objective weight} \\
\hline
\endhead

\hline \multicolumn{3}{r}{\textit{continued on next page}} \\
\endfoot
\endlastfoot %

$\varDutyAssignment_{p, d}$ &
Binary variable equal to 1 if physician~$p$ is assigned to duty~$d$, 0 otherwise &
$c_{p,d}$ \\ \hline

$\varShiftAssignment_{p, s}$ &
Binary variable equal to 1 if physician~$p$ is assigned to shift~$s$, 0 otherwise &
$c_{p,s}$ \\ \hline 

$\varShiftDesiredNumber_{s}$ &
Nonnegative integer variable representing the number of physicians assigned to shift~$s$ above the required minimum~$\parShiftMinimumNumber(s)$, up to the desired minimum~$\parShiftDesiredNumber(s)$ &
$c^{\textnormal{desired}}_s$ \\ \hline

$\varShiftMaximalNumber_{s}$ &
Nonnegative integer variable representing the number of physicians assigned to shift~$s$ above the desired minimum~$\parShiftDesiredNumber(s)$, up to the maximal allowed number~$\parShiftMaxNumber(s)$ &
$c^{\max}_s$ \\ \hline  

$\varShiftDesiredNumberAuxiliary_{s}$ &
Binary auxiliary variable equal to~1 if $\varShiftDesiredNumber_s$ is \emph{not} at its upper bound $\parShiftDesiredNumber(s)$ 
& -- \\ \hline  

$\varDutyBlockAssignment_{p, b}$ &
Binary variable equal to 1 if physician~$p$ is assigned to duty block~$b$, 0 otherwise &
-- \\ \hline  

$\varShiftBlockAssignment_{p, b}$ &
Binary variable equal to 1 if physician~$p$ is assigned to shift block~$b$, 0 otherwise &
-- \\ \hline  

$\varShiftBlockConsecPeriods_{p, b}$ &
Binary variable equal to 1 if physician~$p$ is assigned to shift block~$b$ and the preceding shift block~$\parPrev(b)$ that should be assigned to the same physicians &
$c^{\text{block\_consec}}_{p,b}$ \\ \hline  

$\varDutyConsecPeriods_{p, d}$ &
Binary variable equal to 1 if physician~$p$ is assigned to duty~$d$ and the preceding duty~$\parPrev(d)$ that should be assigned to the same physician &
$c^{\text{consec}}_{p,d}$ \\ \hline  

$\varVioRestTime_{p, d_1, d_2}$ &
Binary variable equal to 1 if physician~$p$ is assigned to both duties/shifts~$d_1$ and~$d_2$, violating a desired rest time &
$c^{\text{rest}}_{p,s1,s2}$ \\ \hline  

$\varWeekendAttendence_{p, w}$ &
Binary variable equal to 1 if physician~$p$ works any duty on weekend~$w$, 0 otherwise &
-- \\ \hline  

$\varWeekendViolation_{p, w}$ &
Nonnegative integer variable counting the degree to which physician~$p$’s weekend preference is violated on weekend~$w$ &
$c^{\text{we\_pref}}_{p}$ \\ \hline  

$\varMaxWeekendViolation_{p, m}$ &
Nonnegative integer variable counting the number of worked weekends of physician~$p$ in month~$m$ above the desired maximum~$\parDesiredMaximumWeekends$ &
$c^{\text{max\_we}}_{p,m}$ \\ \hline  

$\varMinFreeWeekendViolation_{p, m}$ &
Nonnegative integer variable counting the number of free weekends physician~$p$ has fewer than the desired minimum~$\parDesiredMinFreeWeekends$ in month~$m$ &
$c^{\text{free\_we}}_{p,m}$ \\ \hline  

$\varMaxNumberOfDutiesViolation_{p, \pi}$ &
Nonnegative integer variable counting the number of duties in pool~$\pi$ assigned to physician~$p$ above the desired maximum~$\parDesiredMaxNumberOfDutiesInPool(\pi)$ &
$c^{\text{max\_d}}_{\pi,m}$ \\ \hline  

$\varMinNumberOfDutiesViolation_{p, \pi}$ &
Nonnegative integer variable counting the number of duties in pool~$\pi$ assigned to physician~$p$ below the desired minimum~$\parDesiredMinNumberOfDutiesInPool(\pi)$ &
$c^{\text{min\_d}}_{\pi,m}$ \\ \hline  

$\varMaxNumberOfPhysiciansViolation_{\pi, t}$ &
Nonnegative integer variable counting the number of physicians in pool~$\pi$ assigned to duties in pool~$\pi$ on day~$t$ above the allowed maximum~$\parMaxPhysiciansPoolPerDay(\pi)$ &
$c^{\text{max\_phy}}_{\pi,t}$ \\ \hline  

$\varDeviationDutiesPoolNeg_{p, \pi}$ &
Nonnegative integer variable counting the shortfall of duties assigned to physician~$p$ in pool~$\pi$ relative to $\lfloor\parFairDistribution(\pi, p)\rfloor$ &
$c^-_{\pi}$ \\ \hline  

$\varDeviationDutiesPoolPos_{p, \pi}$ &
Nonnegative integer variable counting the excess of duties assigned to physician~$p$ in pool~$\pi$ relative to $\lceil\parFairDistribution(\pi, p)\rceil$ &
$c^+_{\pi}$ \\ \hline  

$\varMaxConsSBViolation_j$ &
Binary variable equal to 1 if the same physician is assigned to all shift blocks in set~$\setConsShiftBlocks_j$, violating the maximum consecutive assignment constraint &
$c^{\text{max\_cons\_b}}_{j}$ \\ \hline 
\end{longtable}

\noindent
\textbf{Objective Function:}
\medskip

\noindent
The linear objective function to be maximized is a weighted sum of all decision variables, with their respective weights (objective coefficients) specified in the column \emph{Objective Weights} in Table~\ref{tab:variables}. If a coefficient has less indices than the corresponding variable, this coefficient is identical for all values of the missing index. As shown in Table~\ref{tab:assignment_coeffs}, the coefficients~$c_{p, d}$ of the duty assignment variables~$\varDutyAssignment_{p, d}$ and $c_{p, s}$ of the shift assignment variables~$\varShiftAssignment_{p, s}$ each consist of different summands that encourage or penalize assigning physician~$p$ to duty~$d$ or shift~$s$ based on different aspects such as desired qualifications or the physician's individual preferences.

\begin{table}[h!]
\centering
\small
\begin{tabularx}{\linewidth}{lX}
\toprule
\textbf{Coefficient} & \textbf{Description} \\
\midrule
$c^1_d$ & Encourages assigning a physician to a non-mandatory duty $d \in \setDuties \setminus \setMandatoryDuties$. \\
$c^2_{p,d}$ & Penalizes the assignment of a physician~$p$ to a duty~$d$ for which she does not have all desired qualifications, i.e., $d \in \setDuties \setminus \setSoftQualifiedDuties(p)$. \\
$c^3_{p,d}$ & Penalizes duty assignments conflicting with soft rest times from the previous planning period, i.e., $d \in \setDutySoftRestLastPP(p)$. \\
$c^4_{p,d}$ & Encourages or penalizes duty assignments based on individual preferences of physician~$p$. \\
\midrule
$c^1_{p,s}$ & Penalizes the assignment of a physician~$p$ to a shift~$s$ for which she does not have all desired qualifications, i.e., $s \in \setShifts \setminus \setSoftQualifiedShift(p)$. \\
$c^2_{p,s}$ & Penalizes shift assignments conflicting with soft rest times from the previous planning period, i.e., $s \in \setShiftSoftRestLastPP(p)$. \\
$c^3_{p,s}$ & Encourages or penalizes shift assignments based on individual preferences of physician~$p$. \\
\bottomrule
\end{tabularx}
\caption{Decomposition of objective coefficients for duty assignment variables~$\varDutyAssignment_{p, d}$ (upper part) and shift assignment variables~$\varShiftAssignment_{p, s}$ (lower part). 
Positive coefficients encourage (i.e., reward) assignments, while negative coefficients penalize them.}
\label{tab:assignment_coeffs}
\end{table}

\medskip

\noindent
\textbf{Constraints:}
\medskip

\noindent
The constraints are formulated verbally, while their mathematical formulation can be found in~\ref{PR:sec:appendix_MIP_formulation}. 

\medskip
\noindent
\textbf{Duty Assignment:} 
Mandatory / non-mandatory duties must be assigned to exactly / at most one physician, respecting existing manual assignments.

\begin{enumerate}[label=(\arabic*), ref=\theenumi]
    \item \label{PR:con:mandatory_duty} \largeStepBackSmall\hyperref[PR:app:mandatory_duty]{\stepForLinkSizeSmall}\stepToGoBackToOriginSmall 
    Each mandatory duty~$d \in \setMandatoryDuties$ must be assigned to exactly one physician~$p \in \setPhysicians$.

    \item \label{PR:con:max_one_duty} \largeStepBackSmall\hyperref[PR:app:max_one_duty]{\stepForLinkSizeSmall}\stepToGoBackToOriginSmall 
    Each non-mandatory duty~$d \in \setDuties \setminus \setMandatoryDuties$ can be assigned to at most one physician~$p \in \setPhysicians$.

    \item \label{PR:con:fixed_duty} \largeStepBackSmall\hyperref[PR:app:fixed_duty]{\stepForLinkSizeSmall}\stepToGoBackToOriginSmall 
    The duties~$d \in \setManuallyPlannedDutyAssignments(p)$ are manually assigned to physician~$p$.

    \item \label{PR:con:manual_duty} \largeStepBackSmall\hyperref[PR:app:manual_duty]{\stepForLinkSizeSmall}\stepToGoBackToOriginSmall 
    A manually planned physician~$p \in \setManuallyPlannedPhysicians$ cannot be assigned to any duty~$d \in \setDuties \setminus \setManuallyPlannedDutyAssignments(p)$ that was not assigned manually.
\end{enumerate}

\medskip
\noindent
\textbf{Qualifications:} 
Physicians can only be assigned to duties and shifts for which they are qualified.

\begin{enumerate}[resume, label=(\arabic*), ref=\theenumi]
    \item \label{PR:con:qualification_duties} \largeStepBack\hyperref[PR:app:qualification_duties]{\stepForLinkSize}\stepToGoBackToOrigin 
    A physician~$p \in \setPhysicians$ cannot be assigned to any duty~$d \in \setDuties \setminus \setQualifiedDuties(p)$ or shift~$s \in \setShifts \setminus \setQualifiedShift(p)$ for which she lacks at least one required qualification or holds any excluded qualification.
\end{enumerate}

\medskip
\noindent
\textbf{Shift Assignment:} 
Each shift must satisfy minimum and maximum staffing levels and adhere to qualification and manual-planning rules.

\begin{enumerate}[resume, label=(\arabic*), ref=\theenumi]
    \item \label{PR:con:shift_1} \largeStepBackSmall\hyperref[PR:app:shift_1]{\stepForLinkSizeSmall}\stepToGoBackToOriginSmall 
    Each shift~$s \in \setShifts$ must have at least $\parShiftMinimumNumber(s)$ physicians assigned.

    \item \label{PR:con:shift_2} \largeStepBackSmall\hyperref[PR:app:shift_2]{\stepForLinkSizeSmall}\stepToGoBackToOriginSmall 
    Each shift~$s \in \setShifts$ can have at most $\parShiftMaxNumber(s)$ physicians assigned.

    \item \label{PR:con:shift_3} \largeStepBackSmall\hyperref[PR:app:shift_3]{\stepForLinkSizeSmall}\stepToGoBackToOriginSmall
    For each shift~$s \in \setShifts$, the sum of $\varShiftDesiredNumber_{s}$ and $\varShiftMaximalNumber_{s}$ equals the number of physicians assigned beyond the required minimum~$\parShiftMinimumNumber(s)$.

    \item \label{PR:con:upper_bound_shift_desired} \largeStepBackSmall\hyperref[PR:app:upper_bound_shift_desired]{\stepForLinkSizeSmall}\stepToGoBackToOriginSmall 
    For each shift~$s \in \setShifts$, the variable~$\varShiftDesiredNumber_{s}$ is upper-bounded by~$\parShiftDesiredNumber(s) - \parShiftMinimumNumber(s)$.

    \item \label{PR:con:upper_bound_shift_max} \largeStepBackSmall\hyperref[PR:app:upper_bound_shift_max]{\stepForLinkSizeSmall}\stepToGoBackToOriginSmall 
    For each shift~$s \in \setShifts$, the variable~$\varShiftMaximalNumber_{s}$ is upper-bounded by~$\parShiftMaxNumber(s) - \parShiftDesiredNumber(s)$.

    \item \label{PR:con:shift_4} \largeStepBack\hyperref[PR:app:shift_4]{\stepForLinkSize}\stepToGoBackToOrigin 
    For each shift~$s \in \setShifts$, the variable~$\varShiftMaximalNumber_{s}$ can take a positive value only if~$\varShiftDesiredNumber_{s}$ is at its upper bound.

    \item \label{PR:con:shift_6} \largeStepBack\hyperref[PR:app:shift_6]{\stepForLinkSize}\stepToGoBackToOrigin 
    Only physicians belonging to the corresponding ward~$\setShiftPhysician(s)$ can be assigned to shift~$s \in \setShifts$.

    \item \label{PR:con:fixed_shift} \largeStepBack\hyperref[PR:app:fixed_shift]{\stepForLinkSize}\stepToGoBackToOrigin 
    The shifts~$s \in \setManuallyPlannedShiftAssignments(p)$ are manually assigned to physician~$p$.

    \item \label{PR:con:man_planned_physicians_exclude} \largeStepBack\hyperref[PR:app:man_planned_physicians_exclude]{\stepForLinkSize}\stepToGoBackToOrigin 
    A manually planned physician~$p \in \setManuallyPlannedPhysicians$ cannot be assigned to any shift~$s \in \setShifts \setminus \setManuallyPlannedShiftAssignments(p)$ that was not assigned manually.
\end{enumerate}

\medskip
\noindent
\textbf{Rest Times:} 
Assignments must / should respect mandatory / desired rest times between consecutive duties or shifts.

\begin{enumerate}[resume, label=(\arabic*), ref=\theenumi]
    \item \label{PR:con:hard_rest_duties} \largeStepBack\hyperref[PR:app:hard_rest_duties]{\stepForLinkSize}\stepToGoBackToOrigin 
    No physician~$p \in \setPhysicians$ can be assigned to a pair of duties or shifts that violates a mandatory rest time, i.e., $(d_1, d_2) \in \setConflicts$ with~$d_1, d_2 \in \setDuties \cup \setShifts$.

    \item \label{PR:con:soft_rest_duties} \largeStepBack\hyperref[PR:app:soft_rest_duties]{\stepForLinkSize}\stepToGoBackToOrigin 
    The variable~$\varVioRestTime_{p, d_1, d_2}$ equals~1 if physician~$p$ is assigned to a pair of duties or shifts that violates a desired rest time, i.e., $(d_1, d_2) \in \setConflictsSoft$ with~$d_1, d_2 \in \setDuties \cup \setShifts$.
\end{enumerate}

\medskip
\noindent
\textbf{Absences:} 
Physicians cannot be assigned to duties or shifts during their absences, or to specific duties immediately before or after them.

\begin{enumerate}[resume, label=(\arabic*), ref=\theenumi]
    \item \label{PR:con:absence_duty} \largeStepBack\hyperref[PR:app:absence_duty]{\stepForLinkSize}\stepToGoBackToOrigin 
    No physician~$p \in \setPhysicians$ is assigned to any duty~$d \in \setDuties(t)$ or shift~$s \in \setShifts(t)$ on an absent day~$t \in \setAbsences(p)$.

    \item \label{PR:con:absence_duty_specific} \largeStepBack\hyperref[PR:app:absence_duty_specific]{\stepForLinkSize}\stepToGoBackToOrigin 
    No physician~$p \in \setPhysicians$ is assigned to any impossible duty~$d \in \setDutyAbsences(p)$ or shift~$s \in \setShiftAbsences(p)$.

    \item \label{PR:con:duties_before_absence} \largeStepBack\hyperref[PR:app:duties_before_absence]{\stepForLinkSize}\stepToGoBackToOrigin 
    No physician~$p \in \setPhysicians$ is assigned to any duty~$d \in \setDutiesBeforeAbsence \cap \setDuties(t-1)$ on the day before an absent day~$t \in \setAbsences(p)$.

    \item \label{PR:con:duties_after_absence} \largeStepBack\hyperref[PR:app:duties_after_absence]{\stepForLinkSize}\stepToGoBackToOrigin 
    No physician~$p \in \setPhysicians$ is assigned to any duty~$d \in \setDutiesAfterAbsence \cap \setDuties(t+1)$ on the day after an absent day~$t \in \setAbsences(p)$.
\end{enumerate}

\medskip
\noindent
\textbf{Duty and Shift Blocks:} 
All duties or shifts within a block must be assigned consistently, and required free days after blocks must be respected.

\begin{enumerate}[resume, label=(\arabic*), ref=\theenumi]
    \item \label{PR:con:duty_block_assignment} \largeStepBack\hyperref[PR:app:duty_block_assignment]{\stepForLinkSize}\stepToGoBackToOrigin 
    All duties~$d \in b$ within a duty block~$b \in \setDutyBlocks$ must be assigned to the same physician.

    \item \label{PR:con:shift_block_assignment} \largeStepBack\hyperref[PR:app:shift_block_assignment]{\stepForLinkSize}\stepToGoBackToOrigin 
    All shifts~$s \in b$ within a shift block~$b \in \setShiftBlocks$ must be assigned to the same physician.

    \item \label{PR:con:free_days_after_duty_block_assignment} \largeStepBack\hyperref[PR:app:free_days_after_duty_block_assignment]{\stepForLinkSize}\stepToGoBackToOrigin 
    Each physician must have at least $\parRequiredFreeDaysAfterDutyBlock(b)$ free days after a duty or shift block~$b \in \setDutyBlocks \cup \setShiftBlocks$.

    \item \label{PR:con:no_other_duties_in_duty_block} \largeStepBack\hyperref[PR:app:no_other_duties_in_duty_block]{\stepForLinkSize}\stepToGoBackToOrigin 
    No physician can be assigned to additional duties during a duty or shift block~$b \in \setBlocksNoAddDuty$, or to additional shifts during a block~$b \in \setBlocksNoAddShift$.
\end{enumerate}

\medskip
\noindent
\textbf{Consecutive Assignments:} 
Assignments across consecutive duties and blocks should comply with continuity and maximum-length requirements.

\begin{enumerate}[resume, label=(\arabic*), ref=\theenumi]
    \item \label{PR:con:shift_block_consec_periods} \largeStepBack\hyperref[PR:app:shift_block_consec_periods]{\stepForLinkSize}\stepToGoBackToOrigin 
    If~$\parPrev(b)$ is not \texttt{None}, the shift block~$b \in \setShiftBlocks$ should be assigned to the same physicians as the preceding block~$\parPrev(b)$.

    \item \label{PR:con:max_consecutive_shift_blocks} \largeStepBack\hyperref[PR:app:max_consecutive_shift_blocks]{\stepForLinkSize}\stepToGoBackToOrigin 
    For each set~$\setConsShiftBlocks_j \in \setConsShiftBlocks$ of consecutive shift blocks that should not all be assigned to the same physician, no physician should be assigned to all blocks in~$\setConsShiftBlocks_j$.

    \item \label{PR:con:consecutive_duties_same_physician} \largeStepBack\hyperref[PR:app:consecutive_duties_same_physician]{\stepForLinkSize}\stepToGoBackToOrigin 
    If~$\parPrev(d)$ is not \texttt{None}, the duty~$d \in \setDuties$ should be assigned to the same physician as the preceding duty~$\parPrev(d)$.
\end{enumerate}

\medskip
\noindent
\textbf{Pools:} 
Duties in pools must respect individual and collective limits and ensure a fair distribution.

\begin{enumerate}[resume, label=(\arabic*), ref=\theenumi]
    \item \label{PR:con:exact_number_duties_pool} \largeStepBack\hyperref[PR:app:exact_number_duties_pool]{\stepForLinkSize}\stepToGoBackToOrigin 
    If~$\parExactDutiesPoolPerMonth(\pi)$ is not \texttt{None}, each physician~$p \in \setPhysicianPool(\pi)$ must be assigned exactly~$\parExactDutiesPoolPerMonth(\pi)$ duties in pool~$\pi \in \setPools$.

    \item \label{PR:con:max_number_duties_pool} \largeStepBack\hyperref[PR:app:max_number_duties_pool]{\stepForLinkSize}\stepToGoBackToOrigin 
    If~$\parMaxNumberOfDutiesInPool(\pi)$ is not \texttt{None}, the total number of duties in pool~$\pi \in \setPools$ assigned to physician~$p \in \setPhysicians$ must not exceed~$\parMaxNumberOfDutiesInPool(\pi)$.

    \item \label{PR:con:desired_max_number_duties_pool} \largeStepBack\hyperref[PR:app:desired_max_number_duties_pool]{\stepForLinkSize}\stepToGoBackToOrigin 
    If~$\parDesiredMaxNumberOfDutiesInPool(\pi)$ is not \texttt{None}, the variable~$\varMaxNumberOfDutiesViolation_{p, \pi}$ equals the number of duties in pool~$\pi \in \setPools$ assigned to physician~$p \in \setPhysicians(\pi)$ above the desired maximum~$\parDesiredMaxNumberOfDutiesInPool(\pi)$.

    \item \label{PR:con:min_number_duties_pool} \largeStepBack\hyperref[PR:app:min_number_duties_pool]{\stepForLinkSize}\stepToGoBackToOrigin 
    If~$\parMinNumberOfDutiesInPool(\pi)$ is not \texttt{None}, the total number of duties in pool~$\pi \in \setPools$ assigned to physician~$p \in \setPhysicians$ must be at least~$\parMinNumberOfDutiesInPool(\pi)$.

    \item \label{PR:con:desired_min_number_duties_pool} \largeStepBack\hyperref[PR:app:desired_min_number_duties_pool]{\stepForLinkSize}\stepToGoBackToOrigin 
    If~$\parDesiredMinNumberOfDutiesInPool(\pi)$ is not \texttt{None}, the variable~$\varMinNumberOfDutiesViolation_{p, \pi}$ equals the number of duties in pool~$\pi$ assigned to physician~$p \in \setPhysicians(\pi)$ below the desired minimum~$\parDesiredMinNumberOfDutiesInPool(\pi)$.

    \item \label{PR:con:max_physicians_pool_day} \largeStepBack\hyperref[PR:app:max_physicians_pool_day]{\stepForLinkSize}\stepToGoBackToOrigin 
    If~$\parMaxPhysiciansPoolPerDay(\pi)$ is not \texttt{None}, the total number of physicians assigned to duties in pool~$\pi \in \setPools$ on any day~$t \in \setDays$ must not exceed~$\parMaxPhysiciansPoolPerDay(\pi)$.

    \item \label{PR:con:desired_max_physicians_pool_day} \largeStepBack\hyperref[PR:app:desired_max_physicians_pool_day]{\stepForLinkSize}\stepToGoBackToOrigin 
    If~$\parMaxPhysiciansPoolPerDay(\pi)$ is not \texttt{None}, the variable~$\varMaxNumberOfPhysiciansViolation_{\pi, t}$ equals the number of physicians in pool~$\pi$ assigned to duties in pool~$\pi \in \setPools$ on day~$t$ above the desired maximum~$\parMaxPhysiciansPoolPerDay(\pi)$.

    \item \label{PR:con:fair_distribution_pool_min} \largeStepBack\hyperref[PR:app:fair_distribution_pool_min]{\stepForLinkSize}\stepToGoBackToOrigin 
    If~$\parFairDistribution(\pi, p)$ is not \texttt{None}, the variable~$\varDeviationDutiesPoolNeg_{p, \pi}$ equals the number of duties in pool~$\pi$ assigned to physician~$p \in \setPhysicians(\pi)$ below the rounded-down target~$\lfloor \parFairDistribution(\pi, p) \rfloor$.

    \item \label{PR:con:fair_distribution_pool_max} \largeStepBack\hyperref[PR:app:fair_distribution_pool_max]{\stepForLinkSize}\stepToGoBackToOrigin 
    If~$\parFairDistribution(\pi, p)$ is not \texttt{None}, the variable~$\varDeviationDutiesPoolPos_{p, \pi}$ equals the number of duties in pool~$\pi$ assigned to physician~$p \in \setPhysicians(\pi)$ above the rounded-up target~$\lceil \parFairDistribution(\pi, p) \rceil$.
\end{enumerate}

\medskip
\noindent
\textbf{Weekends:} 
Weekend assignments must / should respect hard / soft upper and lower bounds, and preferences of physicians for one/ multiple duties per weekend should be respected.

\begin{enumerate}[resume, label=(\arabic*), ref=\theenumi]
    \item \label{PR:con:weekend_attendance_duties} \largeStepBack\hyperref[PR:app:weekend_attendance_duties]{\stepForLinkSize}\stepToGoBackToOrigin 
    The weekend attendance variable~$\varWeekendAttendence_{p, w}$ equals~1 if physician~$p$ is assigned to at least one duty~$d \in \setDutiesOfWeekend(w)$ during weekend~$w \in \setWeekend$.

    \item \label{PR:con:max_consecutive_weekends} \largeStepBack\hyperref[PR:app:max_consecutive_weekends]{\stepForLinkSize}\stepToGoBackToOrigin 
    Each physician~$p \in \setPhysicians$ can work at most~$\parDesiredMaximumConsecutiveWeekends$ consecutive weekends.

    \item \label{PR:con:weekend_preferences_single} \largeStepBack\hyperref[PR:app:weekend_preferences_single]{\stepForLinkSize}\stepToGoBackToOrigin 
    If a physician~$p \in \setPhysicianOneDutyWeekend$ is assigned to multiple duties during weekend~$w \in \setWeekend$, the variable~$\varWeekendViolation_{p, w}$ represents the number of assigned duties minus one.

    \item \label{PR:con:weekend_preferences_multiple} \largeStepBack\hyperref[PR:app:weekend_preferences_multiple]{\stepForLinkSize}\stepToGoBackToOrigin 
    If a physician~$p \in \setPhysicianSeveralDutiesWeekend$ is assigned to only one duty during weekend~$w \in \setWeekend$,
    the variable~$\varWeekendViolation_{p, w}$ equals~1.

    \item \label{PR:con:max_weekends_per_month} \largeStepBack\hyperref[PR:app:max_weekends_per_month]{\stepForLinkSize}\stepToGoBackToOrigin 
    Each physician~$p \in \setPhysicians$ can be assigned to at most $\operatorname{round}(\parMaximumWeekends\cdot\MonthlyWeekendFactor(m))$ weekends in month~$m \in \setMonths$.

    \item \label{PR:con:desired_max_weekends_per_month} \largeStepBack\hyperref[PR:app:desired_max_weekends_per_month]{\stepForLinkSize}\stepToGoBackToOrigin 
    The variable~$\varMaxWeekendViolation_{p, m}$ equals the number of worked weekends of physician~$p \in \setPhysicians$ in month~$m \in \setMonths$ above the desired maximum. 

    \item \label{PR:con:min_free_weekends_per_month} \largeStepBack\hyperref[PR:app:min_free_weekends_per_month]{\stepForLinkSize}\stepToGoBackToOrigin 
    Each physician~$p \in \setPhysicians$ must have at least~$\lfloor\lvert\setWeekendInMonth(m)\lvert - \parDesiredMinFreeWeekends \cdot \MonthlyWeekendFactor(m)\rceil$ free weekends in each month~$m \in \setMonths$.

    \item \label{PR:con:desired_min_free_weekends_per_month} \largeStepBack\hyperref[PR:app:desired_min_free_weekends_per_month]{\stepForLinkSize}\stepToGoBackToOrigin 
    The variable~$\varMinFreeWeekendViolation_{p, m}$ equals the number of free weekends of physician~$p \in \setPhysicians$ in month~$m \in \setMonths$ below the desired minimum. 
\end{enumerate}

\medskip
\noindent
\textbf{Previous Planning Period:} 
Assignments must respect rest times following the previous planning period.

\begin{enumerate}[resume, label=(\arabic*), ref=\theenumi]
    \item \label{PR:con:previous_duty_rest} \largeStepBack\hyperref[PR:app:previous_duty_rest]{\stepForLinkSize}\stepToGoBackToOrigin 
    No physician~$p \in \setPhysicians$ can be assigned to any duty~$d \in \setDutyRestLastPP(p)$ or shift~$s \in \setShiftRestLastPP(p)$ that violates mandatory rest times after duties or shifts from the previous planning period.
\end{enumerate}
\section{Practical Implementation}\label{sec:implementation}
Our model is implemented as a full-featured web application based on the open-source Python framework Flask~\cite{Flask}. 
This architecture enables seamless interaction between a web frontend, providing an intuitive and user-friendly graphical user interface (GUI), and a Python backend responsible for data management and optimization. 
All data are stored securely in an SQLite database. 
The software is deployed on a web server, providing easy access for hospital staff through personal logins for all physicians and other personnel (e.g., department secretaries) involved in the rostering process.

Users are divided into two main groups: \emph{physicians}, who can specify their personal preferences and view published rosters to see their assigned duties and shifts; and \emph{planners}, who can additionally define hospital-specific planning requirements and parameters. 
Importantly, planners are not required to have any expertise in mathematical modeling. 
The GUI of our web application allows them to specify every detail of the department-specific roster structure, while the corresponding MIP model is automatically generated and solved in the backend. 
Consequently, the entire roster generation process is guided by the application—from defining the roster structure and physicians’ individual preferences to publishing the final roster, making manual adjustments (e.g., to account for unplanned absences after publication), and even automatically importing duties and shifts into personal calendars.

The tool is highly adaptable and capable of representing the diverse duty and shift structures across the three partner hospitals. 
Furthermore, we have verified with a fourth department specializing in gastroenterology and endocrinology that its rostering problem can also be modeled directly using our software. 
The only anticipated adjustments concern minor GUI customizations desired by the department’s planners; no changes to the underlying MIP model are required.

\begin{figure}[H]
   \centering
   \includegraphics[width=0.4\textwidth]{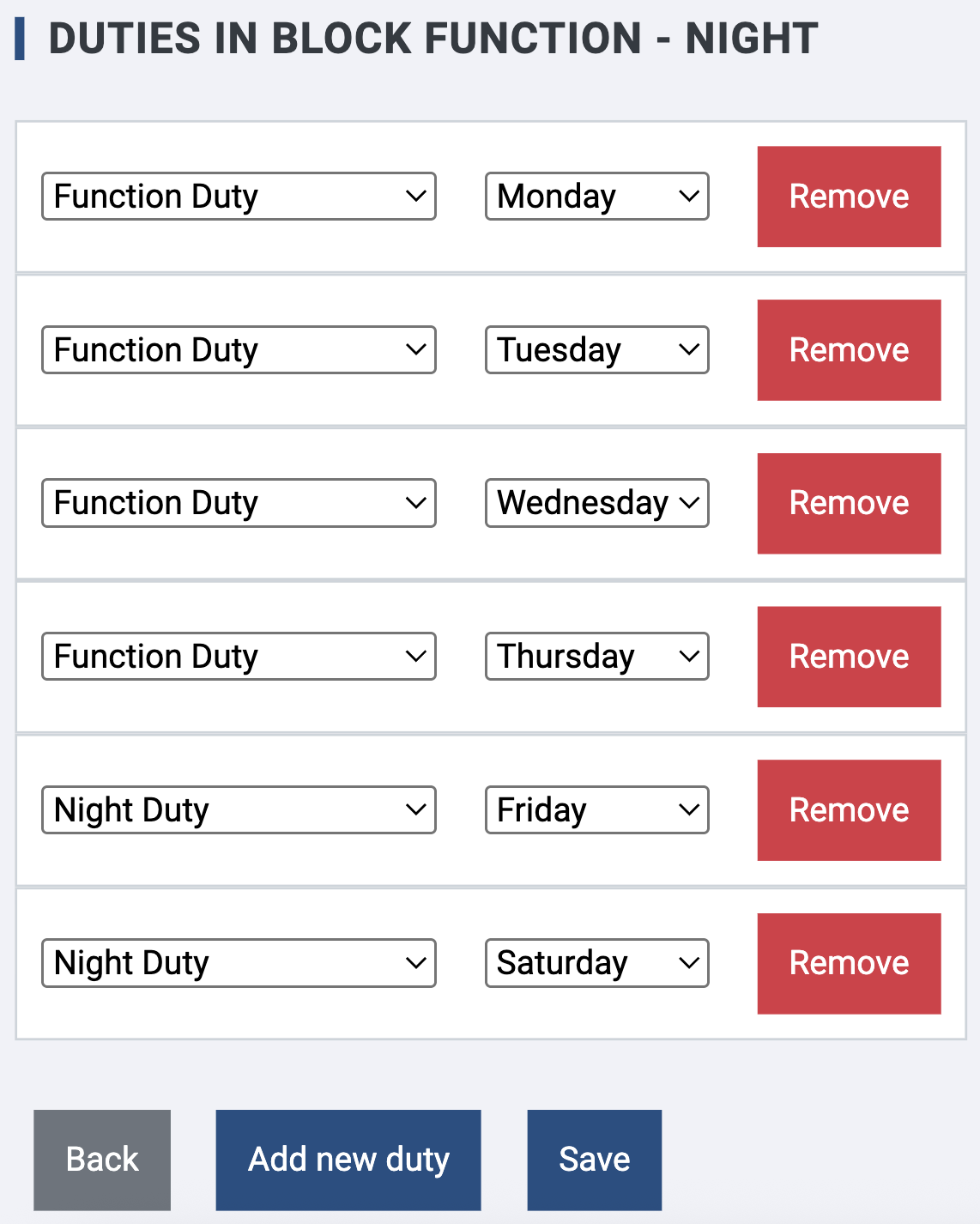}
   \caption{Creation of a duty block consisting of different duties over a week in the web application.}
   \label{PR:fig:block_creation}
\end{figure}

\begin{figure}[H]
   \centering
   \begin{subfigure}[b]{0.48\textwidth}
       \centering
       \includegraphics[width=\textwidth]{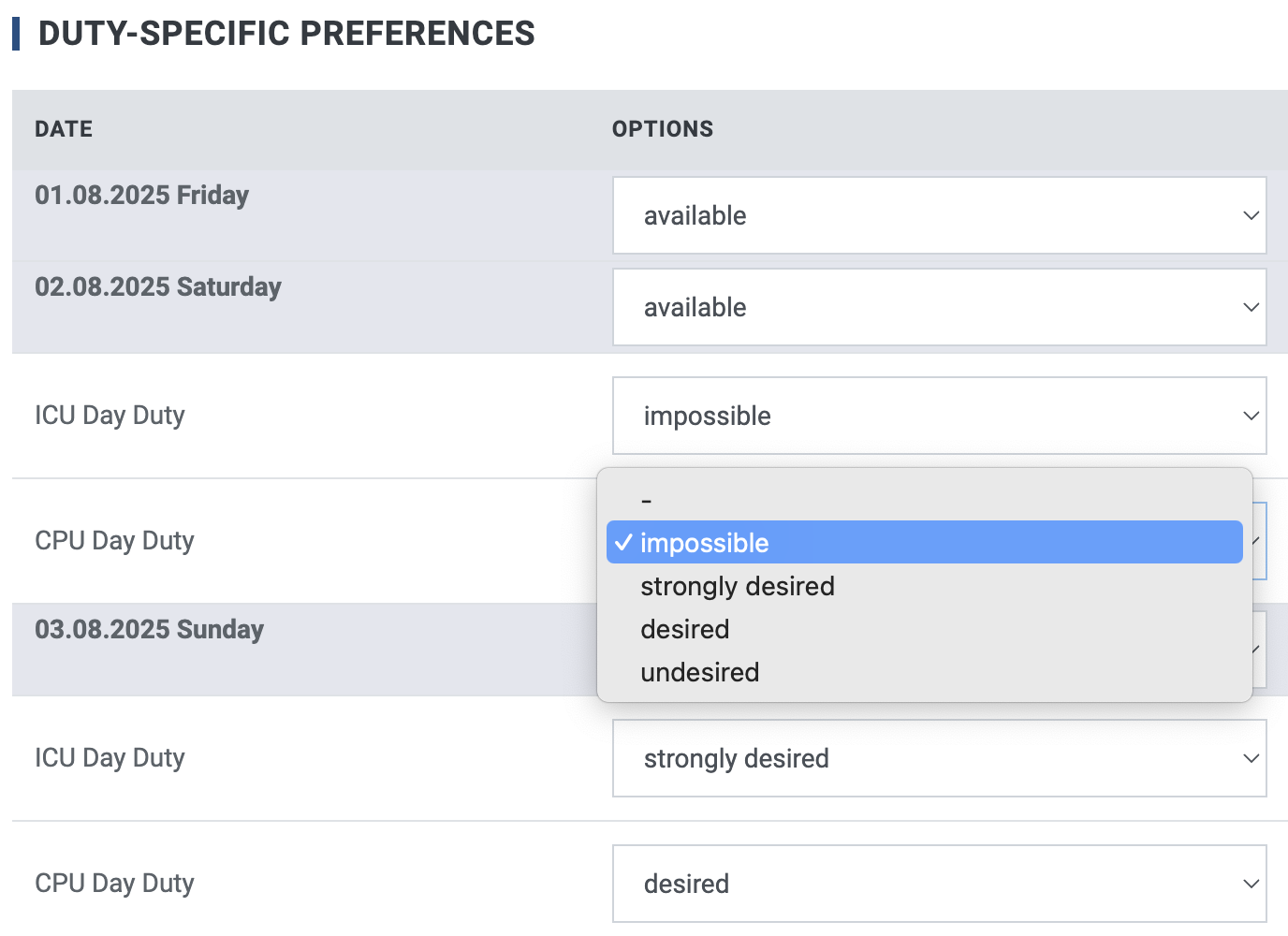}
       \caption{Layout for entering duty-specific preferences.}
       \label{PR:fig:duty_specific_preferences}
   \end{subfigure}
   \hfill
   \begin{subfigure}[b]{0.48\textwidth}
       \centering
       \includegraphics[width=\textwidth]{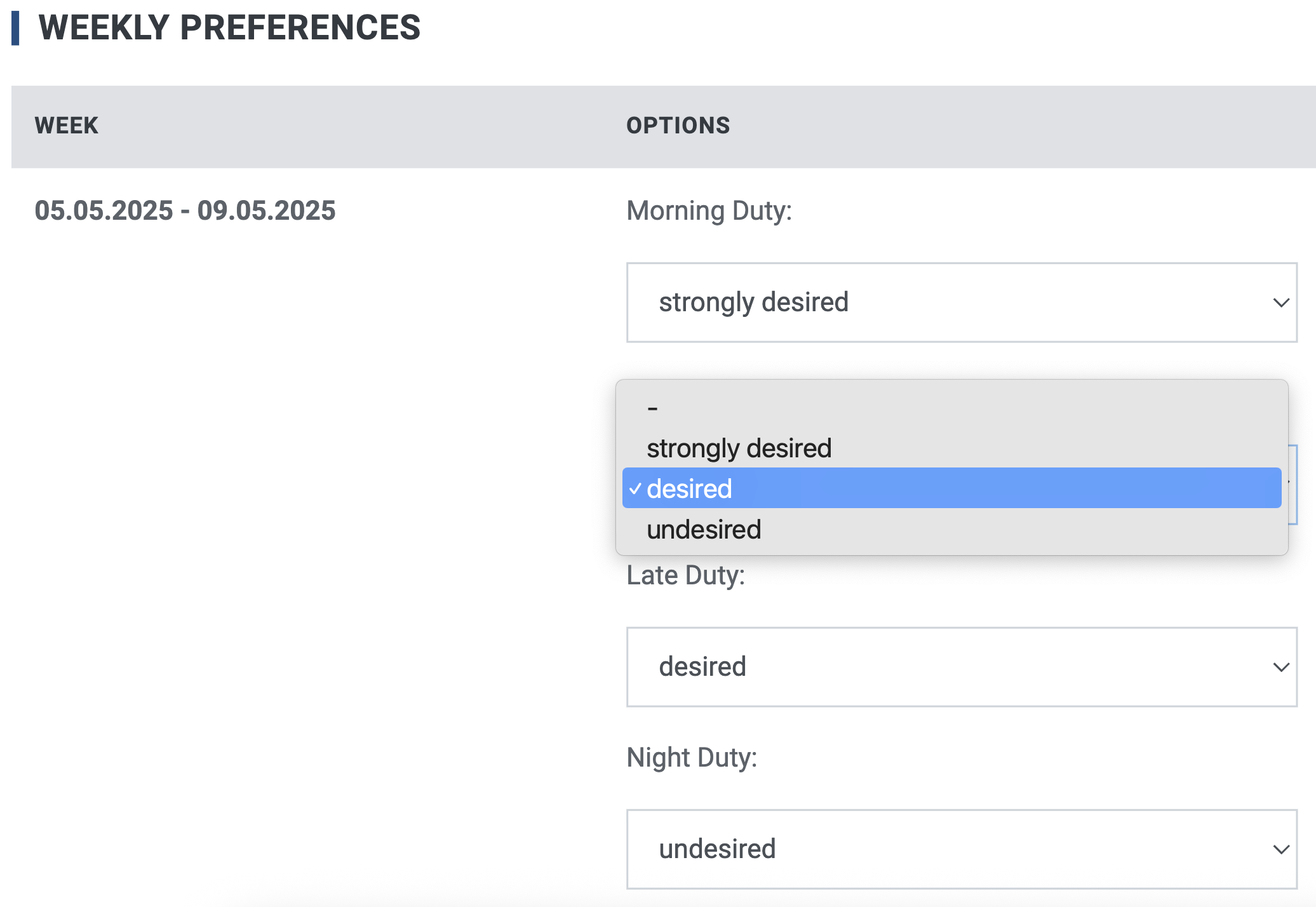}
       \caption{Layout for entering weekly preferences.}
       \label{PR:fig:block_preferences}
   \end{subfigure}
   \caption{Input interfaces for specifying physicians’ preferences in the web application. 
   (a) Duty-specific preferences; (b) Weekly preferences.}
   \label{PR:fig:physicians_preferences}
\end{figure}

To illustrate the GUI, Figure~\ref{PR:fig:block_creation} shows how a duty block containing different duties over a week can be created, while Figure~\ref{PR:fig:physicians_preferences} displays the interface for entering physicians’ preferences. 
The left-hand side shows the layout for duty-specific preferences, whereas the right-hand side presents the version for weekly preferences.

In summary, its implementation as a full-featured web application is a key strength of our approach, which enhances its practical applicability by making the model accessible to users without a mathematical or technical background.  
Moreover, beyond these usability benefits, the architecture offers improved data protection compared to common practices in hospitals, where physicians typically submit their preferences via email or on a physical list posted on-site.

\section{Model Evaluation with Real-World Data}\label{PR:sec:comp_results}

In this section, we evaluate the roster quality and computation times of our general model using four monthly real-world instances (March--June) from each of the three hospital departments---internal medicine, cardiology, and orthopedics/trauma surgery---whose rostering problems are described in Section~\ref{sec:problem_description}.

\medskip

For the cardiology department, which already uses our web application, the data were obtained directly from the application. For the internal me\-di\-cine and orthopedics/trauma surgery departments, data were collected from Excel files used as inputs for department-specific MIP models that have been in use for several years.
For the internal medicine and orthopedics/trauma surgery departments, we use these department-specific models as benchmarks to evaluate the quality of the rosters generated by our general model. 
Quality indicators include violations of staffing requirements, fulfillment of physician preferences, fair distribution of duties, consecutive assignment of single duties and duty/shift blocks (where applicable), and compliance with desired rest times (where applicable). 
For the cardiology department, no benchmark from a department-specific model is available. Nevertheless, the results provide valuable insights, as this department faces significant personnel shortages and frequently requires support from other departments or has unassigned duties.

\medskip

All computational results were obtained on a laptop running macOS Sequoia~15.5, equipped with an Apple~M1 Pro processor (eight cores: six performance cores at 3.2~GHz and two efficiency cores at 2.0~GHz) and 32~GB of memory. MIP models were solved using the open-source solver CBC~2.10.11~\cite{CBC} with an optimality gap of 3\%, and the overall code was executed in Python~3.9.12.

\subsection{Departments with Established Benchmark Models}

We present the evaluation for the internal medicine and orthopedics/trauma surgery departments, comparing our general MIP model with the existing department-specific MIP models that have been in use for several years. The results for the internal medicine department are shown in Table~\ref{PR:tab:analysis_hospital1}, and those for the orthopedics/trauma surgery department in Table~\ref{PR:tab:analysis_hospital3}.

{\small
\begin{longtable}{|l|rrrr>{\cellcolor{gray!20}}r|rrrr>{\cellcolor{gray!20}}r|}
\hline
\multicolumn{11}{|c|}{\textbf{Internal Medicine Department}} \\
\hline
& \multicolumn{5}{c|}{\textbf{Specific Model}} & \multicolumn{5}{c|}{\textbf{General Model}} \\
\hline
& \rotatebox{90}{March} & \rotatebox{90}{April} & \rotatebox{90}{May} & \rotatebox{90}{June} & \rotatebox{90}{Average\;} 
& \rotatebox{90}{March} & \rotatebox{90}{April} & \rotatebox{90}{May} & \rotatebox{90}{June} & \rotatebox{90}{Average\;} \\
\hline
\endfirsthead

\multicolumn{11}{c}%
{\tablename\ \thetable\ -- \textit{continued from previous page}}\\
\hline
& \multicolumn{5}{c|}{\textbf{Specific Model}} & \multicolumn{5}{c|}{\textbf{General Model}} \\
\hline
& \rotatebox{90}{March} & \rotatebox{90}{April} & \rotatebox{90}{May} & \rotatebox{90}{June} & \rotatebox{90}{Average\;} 
& \rotatebox{90}{March} & \rotatebox{90}{April} & \rotatebox{90}{May} & \rotatebox{90}{June} & \rotatebox{90}{Average\;} \\
\hline
\endhead

\hline \multicolumn{11}{r}{\textit{continued on next page}} \\
\endfoot
\endlastfoot %

\multicolumn{11}{|c|}{\textbf{Computation times}} \\
\hline
Solver time [s] & 1 & 1 & 1 & 1 & \textbf{1} & 91 & 66 & 80 & 52 & 72.25 \\
\hline
Total computation time [s] & 3 & 3 & 3 & 3 & \textbf{3} & 145 & 119 & 134 & 108 & 126.5 \\
\hline
\multicolumn{11}{|c|}{\textbf{Staffing requirements}} \\
\hline
\makecell[l]{Understaffed wards} & 9 & 4 & 4 & 8 & 6.25 & 8 & 4 & 5 & 8 & 6.25 \\
\hline
\multicolumn{11}{|c|}{\textbf{Physician preferences}} \\
\hline
\makecell[l]{Number of desired duties} & 75 & 59 & 62 & 57 & 63.25 & 75 & 59 & 62 & 57 & 63.25 \\
\hline
Desired duties assigned & 46 & 33 & 21 & 28 & 32 & 45 & 33 & 23 & 29 & \textbf{32.5} \\
\hline
\makecell[l]{Number of strongly\\ desired duties} & 0 & 26 & 43 & 27 & 24 & 0 & 26 & 43 & 27 & 24 \\
\hline
\makecell[l]{Strongly desired\\ duties assigned} & 0 & 15 & 21 & 15 & \textbf{12.75} & 0 & 15 & 19 & 14 & 12 \\
\hline
Undesired duties assigned & 0 & 1 & 0 & 1 & 0.5 & 0 & 1 & 0 & 1 & 0.5 \\
\hline
\multicolumn{11}{|c|}{\textbf{Fair distribution of duties}} \\
\hline
\makecell[l]{Phys. who want fewer duties} & 13 & 15 & 15 & 11 & 13.5 & 13 & 15 & 15 & 11 & 13.5 \\
\hline
\makecell[l]{Of these, phys. who\\ get fewer duties} & 7 & 8 & 4 & 7 & 6.5 & 10 & 9 & 8 & 9 & \textbf{9} \\
\hline
\makecell[l]{Phys. who want more duties} & 3 & 1 & 2 & 1 & 1.75 & 3 & 1 & 2 & 1 & 1.75 \\
\hline
\makecell[l]{Of these, physicians who\\ get more duties} & 1 & 1 & 1 & 0 & 0.75 & 3 & 1 & 1 & 1 & \textbf{1.5} \\
\hline
\makecell[l]{Phys. who want neither\\ fewer nor more duties} & 4 & 4 & 4 & 3 & 3.75 & 4 & 4 & 4 & 3 & 3.75 \\
\hline
\makecell[l]{Of these, phys. who get\\ neither fewer nor more duties} & 4 & 3 & 1 & 3 & 2.75 & 4 & 2 & 2 & 3 & 2.75 \\
\hline
\makecell[l]{Phys. with too many duties\\ on weekends and holidays} & 11 & 9 & 3 & 12 & 8.75 & 11 & 9 & 3 & 12 & 8.75 \\
\hline
\multicolumn{11}{|c|}{\textbf{Desired rest times}} \\
\hline
\makecell[l]{Assigned duties with only \\ one duty-free day in between} & 0 & 2 & 4 & 1 & 1.75 & 2 & 0 & 1 & 3 & \textbf{1.5} \\
\hline
\makecell[l]{Assigned duties with only \\ two duty-free days in between} & 3 & 2 & 3 & 4 & 3 & 2 & 2 & 4 & 2 & \textbf{2.5} \\
\hline
\caption{Comparison of general and department-specific models for the internal medicine department. Bold values denote strictly better average performance for the corresponding indicator.} \label{PR:tab:analysis_hospital1}
\end{longtable}
}

{\small
\renewcommand{\arraystretch}{1.5}
\begin{longtable}{|l|rrrr>{\columncolor{gray!20}}r|rrrr>{\columncolor{gray!20}}r|}
\hline
\multicolumn{11}{|c|}{\textbf{Orthopedics/Trauma Surgery Department}} \\
\hline
& \multicolumn{5}{c|}{\textbf{Specific Model}} & \multicolumn{5}{c|}{\textbf{General Model}} \\
\hline
& \rotatebox{90}{March} & \rotatebox{90}{April} & \rotatebox{90}{May} & \rotatebox{90}{June} & \rotatebox{90}{Average\;} 
& \rotatebox{90}{March} & \rotatebox{90}{April} & \rotatebox{90}{May} & \rotatebox{90}{June} & \rotatebox{90}{Average\;} \\
\hline
\endfirsthead

\multicolumn{11}{c}%
{\tablename\ \thetable\ -- \textit{continued from previous page}}\\
\hline
& \multicolumn{5}{c|}{\textbf{Specific Model}} & \multicolumn{5}{c|}{\textbf{General Model}} \\
\hline
& \rotatebox{90}{March} & \rotatebox{90}{April} & \rotatebox{90}{May} & \rotatebox{90}{June} & \rotatebox{90}{Average\;} 
& \rotatebox{90}{March} & \rotatebox{90}{April} & \rotatebox{90}{May} & \rotatebox{90}{June} & \rotatebox{90}{Average\;} \\
\hline
\endhead

\hline \multicolumn{11}{r}{\textit{continued on next page}} \\
\endfoot
\endlastfoot %

\multicolumn{11}{|c|}{\textbf{Computation times}} \\
\hline
Solver time [s] & 9 & 14 & 7 & 10 & \textbf{10} & 325 & 155 & 165 & 304 & 237.25 \\
\hline
Total computation time [s] & 17 & 19 & 14 & 18 & \textbf{17} & 502 & 304 & 349 & 496 & 412.75 \\
\hline
\multicolumn{11}{|c|}{\textbf{Staffing requirements}} \\
\hline
\makecell[l]{Understaffed wards} & 0 & 0 & 0 & 7 & 1.75 & 0 & 0 & 0 & 0 & \textbf{0} \\
\hline
\makecell[l]{Understaffed ward fellows} & 2 & 0 & 2 & 12 & \textbf{4} & 1 & 0 & 3 & 14 & 4.5 \\
\hline
\multicolumn{11}{|c|}{\textbf{Physician preferences}} \\
\hline
\makecell[l]{Number of desired duties} & 77 & 68 & 89 & 83 & 79.25 & 77 & 68 & 89 & 83 & 79.25 \\
\hline
\makecell[l]{Desired duties assigned} & 50 & 40 & 40 & 35 & 41.25 & 48 & 46 & 44 & 41 & \textbf{44.75} \\
\hline
Undesired duties assigned & 0 & 3 & 2 & 6 & 2.75 & 2 & 1 & 3 & 5 & 2.75 \\
\hline
\multicolumn{11}{|c|}{\textbf{Fair distribution of duties}} \\
\hline
\makecell[l]{Phys. with fewer duties \\ than desired on \emph{work days}} & 0 & 1 & 0 & 1 & 0.5 & 1 & 0 & 1 & 0 & 0.5 \\
\hline
\makecell[l]{Phys. with more duties \\ than desired on \emph{work days}} & 0 & 0 & 0 & 0 & \textbf{0} & 0 & 1 & 0 & 0 & 0.25 \\
\hline
\makecell[l]{Phys. with fewer duties \\ than desired on \emph{weekends}} & 0 & 0 & 0 & 0 & 0 & 0 & 0 & 0 & 0 & 0 \\
\hline
\makecell[l]{Phys. with more duties \\ than desired on \emph{weekends}} & 0 & 0 & 0 & 0 & 0 & 0 & 0 & 0 & 0 & 0 \\
\hline
\multicolumn{11}{|c|}{\textbf{Consecutive assignments}} \\
\hline
\makecell[l]{Consecutive assignments \\ of late duties} & 24 & 15 & 12 & 14 & 16.25 & 25 & 17 & 10 & 14 & \textbf{16.5} \\
\hline
\caption{Comparison of general and department-specific models for the orthopedics/trauma surgery department. Bold values denote strictly better average performance for the corresponding indicator.}
\label{PR:tab:analysis_hospital3}
\end{longtable}
}

For each instance, the solver time is reported, denoting the duration required by CBC to obtain a solution to the MIP with a 3\% optimality gap. In addition, the total computation time is provided, defined as the elapsed time from the initiation of the computation---either through command-line execution of the model or by triggering the corresponding function in the web application---until the roster generated by each model is completed.

\medskip

We observe that the quality of the rosters produced by the general model is largely identical to that of the department-specific models. While more quality indicators favor the general model overall, no single indicator reveals a significant difference in roster quality.

Regarding computation times, the department-specific models naturally perform faster than the general model. This is due to two factors. First, the specific MIPs are more customized and contain fewer variables (specific model: $\approx$~20,000 vs. general model: $\approx$~400,000) and constraints (specific model: $\approx$~20,000 vs. general model: $\approx$~1,000,000). Second, the creation of the MIPs is significantly faster in these cases, as reading and writing to a database introduces additional overhead compared to using simple Excel files. Nevertheless, the total computation time for the general model remains consistently below nine minutes, and for the internal medicine department, it is always below three minutes. Since physician rostering is not time-critical and occurs only once per month, these computation times are more than satisfactory for practical application.

\enlargethispage{\baselineskip}
Overall, we conclude that the advantages of our model---such as its easy configuration via the GUI and its adaptability to changing roster structures or fairness concepts---far outweigh the slightly increased computation times in practice.

\subsection{Department without a Benchmark Model}

We now present the evaluation for the cardiology department, for which no department-specific benchmark model is available. The results obtained using our general model are shown in Table~\ref{PR:tab:analysis_hospital2}.

Regarding computation times, every instance for the cardiology department is solved in under three minutes of total computation time. When comparing the values of the various quality indicators with those of the other two departments, we generally observe less favorable results. This can be attributed to the department’s severe personnel shortages, which often leave too few physicians available to cover all duties and shifts. In our framework, this is reflected by modeling all duties as \emph{optionally unassigned}, with some remaining unassigned in the generated rosters, while certain wards are understaffed. Moreover, the presence of numerous duty blocks further complicates the planning process by reducing overall scheduling flexibility. Since maximizing duty coverage is prioritized over satisfying individual physician preferences, these shortages also lead to a lower fulfillment rate of personal preferences compared to the other departments.

Nevertheless, our general model succeeds in producing high-quality rosters even under these challenging conditions and has been successfully used in practice by the department for several months. By minimizing the number of unassigned duties and understaffed wards, the model not only helps to limit the need for costly external staff but also contributes to ensuring adequate medical coverage across wards.

{\small
\begin{longtable}{|p{0.5\linewidth}|r|r|r|r|>{\columncolor{gray!20}}r|}
\hline
\multicolumn{6}{|c|}{\textbf{Cardiology Department}} \\
\hline
 & \rotatebox{90}{March} & \rotatebox{90}{April} & \rotatebox{90}{May} & \rotatebox{90}{June} & \rotatebox{90}{Average\;} \\
\hline
\endfirsthead

\multicolumn{6}{c}%
{\tablename\ \thetable\ -- \textit{continued from previous page}} \\
\hline
 & \rotatebox{90}{March} & \rotatebox{90}{April} & \rotatebox{90}{May} & \rotatebox{90}{June} & \rotatebox{90}{Average\;} \\
\hline
\endhead

\hline \multicolumn{6}{r}{\textit{continued on next page}} \\
\endfoot

\hline

\caption{Evaluation results for the cardiology department.} \label{PR:tab:analysis_hospital2}
\endlastfoot

\multicolumn{6}{|c|}{\textbf{Computation times}}\\ \hline
Solver time [s] & 28 & 26 & 29 & 40 & 30.75 \\ \hline
Total computation time [s] & 112 & 112 & 119 & 137 & 120 \\ \hline

\multicolumn{6}{|c|}{\textbf{Staffing requirements}} \\ \hline
Unassigned duties & 13 & 4 & 18 & 26 & 15.25 \\ \hline
Understaffed wards & 0 & 5 & 0 & 4 & 2.25 \\ \hline

\multicolumn{6}{|c|}{\textbf{Physician preferences (duty-specific)}} \\ \hline
Number of desired duties & 24 & 10 & 15 & 15 & 16 \\ \hline
Desired duties assigned & 6 & 2 & 3 & 1 & 3 \\ \hline
Number of strongly desired duties & 39 & 21 & 18 & 15 & 23.25 \\ \hline
Strongly desired duties assigned & 10 & 9 & 8 & 9 & 9 \\ \hline
Undesired duties assigned & 0 & 1 & 2 & 1 & 1 \\ \hline

\multicolumn{6}{|c|}{\textbf{Physician preferences (weekly)}} \\ \hline
Strongly desired weekly sets assigned & 5 & 2 & 8 & 7 & 5.5 \\ \hline
Desired weekly sets assigned & 1 & 1 & 1 & 1 & 1 \\ \hline
Undesired weekly sets assigned & 0 & 1 & 0 & 2 & 0.75 \\ \hline

\multicolumn{6}{|c|}{\textbf{Fair distribution of duties}} \\ \hline
\makecell[l]{Phys. with too few weekend duties} & 3 & 2 & 4 & 8 & 4.25 \\ \hline
\makecell[l]{Phys. with too many weekend duties} & 1 & 1 & 2 & 6 & 2.5 \\ \hline
\makecell[l]{Phys. with more than one \\ ICU night duty block} & 0 & 0 & 0 & 0 & 0 \\ \hline

\multicolumn{6}{|c|}{\textbf{Consecutive assignments and desired rest times}} \\ \hline
\makecell[l]{Consecutive assignments of \\ shift blocks} & 15 & 2 & 12 & 7 & 9 \\ \hline
Violated desired rest times & 0 & 1 & 0 & 1 & 0.5 \\ \hline

\end{longtable}
}

\section{Conclusion and Outlook}\label{PR:sec:conclusion}
Physician rostering is a complex and critical task in hospitals, with fair and high-quality schedules being essential for both operational continuity and physician satisfaction. Many hospitals still rely on manual roster creation, which is time-consuming, error-prone, and often results in unfair schedules or violations of labor regulations.
In contrast to most models in the literature, which are tailored to a single department, we present a general and adaptable framework for physician rostering. Developed and validated with three hospital departments that differ substantially in structure, objectives, and constraints, the model produces high-quality rosters efficiently---even with open-source solvers---and matches or slightly exceeds the quality of established department-specific benchmark models.

The framework is readily adoptable by new departments, and existing users can accommodate changes in planning processes without modifying the underlying optimization model. Preliminary evaluations with additional departments suggest that their rostering requirements can also be represented within the framework.
A key feature is the integration into a web-based application with a graphical user interface, allowing planners to configure parameters flexibly and physicians to view rosters online, ensuring transparency and engagement. One partner department already uses the system in practice, with others expected to follow.

Future research could consider onboarding additional departments or supporting additional roster structures such as cyclic rosters. Further studies could also investigate how physician rostering tools can be more closely linked to hospital information systems to streamline administrative workflows, and how usability can be further improved to increase practical adoption and physician satisfaction.

\section*{Statements and Declarations}

\subsection*{Author contributions}
{\parindent0pt
\textbf{Florian Meier:} Methodology, Software, Validation, Formal analysis, Investigation, Data Curation, Writing - Original Draft, Writing - Review \& Editing, Visualization.
\textbf{Jan Boeckmann:} Methodology, Software, Validation, Formal analysis, Writing - Review \& Editing, Visualization.
\textbf{Clemens Thielen:} Conceptualization, Methodology, Validation, Formal analysis, Writing - Review \& Editing, Supervision, Project Administration, Funding acquisition.
}

\subsection*{Funding}
\noindent
This work was supported by the German Federal Ministry of Research, Technology and Space (BMFTR) -- Grant number 03DPS1158.

\subsection*{Conflicts of interest}
\noindent
The authors have no conflicts of interest to declare that are relevant to the content of this article.

\subsection*{Data availability}
\noindent
\sloppy
The implementation of the mixed-integer programming model, together with the database files containing the (anonymized) real-world instances used for model evaluation, are publicly available at \url{https://github.com/florian995/Physician-Rostering-Framework}.
\fussy
\bibliographystyle{elsarticle-num}
\bibliography{bibfile}

@STRING{EJOR   = "European Journal of Operational Research"}

@STRING{ORHC    = "Operations Research for Health Care"}

@Manual{Flask,
  title        = {Flask: A Python Microframework for Web Development},
  author       = {{A.\ Ronacher}},
  organization = {Pallets Projects},
  year         = {2025},
  version      = {3.1.2},
  note         = {\url{https://pypi.org/project/Flask/}},
}

@Manual{CBC,
  title        = {CBC},
  organization = {COIN-OR Foundation},
  year         = {2025},
  version      = {3.2.2},
  note         = {\url{https://github.com/coin-or/Cbc}},
}

@Article{Rousseau+etal:general-approach,
  author =       {L.-M. Rousseau and G. Pesant and M. Gendreau},
  title =        {A General Approach to the Physician Rostering Problem},
  journal = 	 	 {Annals of Operations Research},
  year = 	 			 {2002},
  volume =	 		 {115},
  pages =	 			 {193--205}
}

@Article{schoenfelder+pfefferlen:german-hospital,
  author =       {J. Schoenfelder and C. Pfefferlen},
  title =        {Decision Support for the Physician Scheduling Process at a {German} Hospital},
  journal = 	 	 {Service Science},
  year = 	 			 {2018},
  volume =	 		 {10},
  number =			 {3},
  pages =	 			 {215--229}
}

@Article{Bruni+Detti:physician-scheduling,
  author =       {R. Bruni and P. Detti},
  title =        {A Flexible Discrete Optimization Approach to the Physician Scheduling Problem},
  journal = 	 	 ORHC,
  year = 	 			 {2014},
  volume =	 		 {3},
  number =	 		 {4},
  pages =	 			 {191--199}
}

@Article{fuchs2025fairness,
  author =       {G. Fuchs and K. Schimmelpfeng and J. O. Brunner},
  title =        {A roadmap for integrating fairness in personnel planning and scheduling in hospitals},
  journal = 	 	 {Operations Research, Data Analytics and Logistics},
  year = 	 			 {2025},
  volume =	 		 {45},
  pages =	 			 {200479}
}

@Article{santos2014,
  author =       {M. Santos and H. Eriksson},
  title =        {Insights into physician scheduling: a case study of public hospital departments in {Sweden}},
  journal = 	 	 {International Journal of Health Care Quality Assurance},
  year = 	 			 {2014},
  volume =	 		 {27},
  number =	 		 {2},
  pages =	 			 {76--90}
}

@Article{cheang2003nurse,
  author =       {B. Cheang and H. Li and A. Lim and B. Rodrigues},
  title =        {Nurse Rostering Problems - A Bibliographic Survey},
  journal = 	 	 EJOR,
  year = 	 			 {2003},
  volume =	 		 {151},
  number =			 {3},
  pages =	 			 {447--460}
}

@Article{Brunner+etal:flexible-shift,
  author =       {J. O. Brunner and J. F. Bard and R. Kolisch},
  title =        {Flexible Shift Scheduling of Physicians},
  journal = 	 	 {Health Care Management Science},
  year = 	 			 {2009},
  volume =	 		 {12},
  number =	 		 {3},
  pages =	 			 {285--305}
}

@article{kraul2024optimizing,
  title={Optimizing physician schedules with resilient break assignments},
  author={Kraul, Sebastian and Erhard, Melanie and Brunner, Jens O},
  journal={Omega},
  volume={129},
  pages={103154},
  year={2024}
}

@article{fugener2019planning,
  title={Planning for overtime: The value of shift extensions in physician scheduling},
  author={F{\"u}gener, Andreas and Brunner, Jens O},
  journal={INFORMS Journal on Computing},
  volume={31},
  number={4},
  pages={732--744},
  year={2019}
}

@article{adams2019physician,
  title={Physician rostering for workload balance},
  author={Adams, Thomas and O’Sullivan, Michael and Walker, Cameron},
  journal={Operations Research for Health Care},
  volume={20},
  pages={1--10},
  year={2019}
}

@article{benazzouz2015literature,
  title={A literature review on the nurses’ planning problems},
  author={Benazzouz, Touria and Echchatbi, Abdelwahed and Bellabdaoui, Adil},
  journal={International Journal of Mathematics and Computational Science},
  volume={1},
  number={5},
  pages={268--274},
  year={2015}
}

@Article{burke2004state,
  author =       {E. K. Burke and P. {De Causmaecker} and G. {Vanden Berghe} and H. {Van Landeghem}},
  title =        {The State of the Art of Nurse Rostering},
  journal = 	 	 {Journal of Scheduling},
  year = 	 			 {2004},
  volume =	 		 {7},
  number =			 {6},
  pages =	 			 {441--499}
}

@article{cappanera2021emergency,
  title={The emergency department physician rostering problem: obtaining equitable solutions via network optimization},
  author={Cappanera, Paola and Visintin, Filippo and Rossi, Roberta},
  journal={Flexible Services and Manufacturing Journal},
  pages={1--44},
  year={2021}
}

@article{wickert2021integer,
  title={An integer programming approach for the physician rostering problem},
  author={Wickert, Toni I and Kummer Neto, Alberto F and Boniatti, M{\'a}rcio M and Buriol, Luciana S},
  journal={Annals of Operations Research},
  volume={302},
  pages={363--390},
  year={2021}
}

@article{abdullah2025advancing,
  title={Advancing Anaesthetist Rostering Quality: {A} Practical Approach Towards Fairness and Efficiency},
  author={Abdullah, Norizal and Ayob, Masri and Lam, Meng Chun and Sabar, Nasser R and Kendall, Graham and Yong, Liu Chian},
  journal={IEEE Access},
  volume={13},
  pages={12692--12708},
  year={2025}
}

@article{dewa2017relationship,
  title={The relationship between physician burnout and quality of healthcare in terms of safety and acceptability: a systematic review},
  author={Dewa, Carolyn S and Loong, Desmond and Bonato, Sarah and Trojanowski, Lucy},
  journal={BMJ Open},
  volume={7},
  number={6},
  pages={e015141},
  year={2017}
}

@article{bowers2016neonatal,
  title={Neonatal physician scheduling at the {University of Tennessee Medical Center}},
  author={Bowers, Melissa R and Noon, Charles E and Wu, Wei and Bass, J Kirk},
  journal={Interfaces},
  volume={46},
  number={2},
  pages={168--182},
  year={2016}
}

@article{thielen2018duty,
  title={Duty rostering for physicians at a department of orthopedics and trauma surgery},
  author={C. Thielen},
  journal={Operations Research for Health Care},
  volume={19},
  pages={80--91},
  year={2018}
}

@article{gross2018online,
  title={Online rescheduling of physicians in hospitals},
  author={Gross, Christopher N and F{\"u}gener, Andreas and Brunner, Jens O},
  journal={Flexible Services and Manufacturing Journal},
  volume={30},
  number={1},
  pages={296--328},
  year={2018}
}

@article{fugener2015duty,
  title={Duty and workstation rostering considering preferences and fairness: a case study at a department of anaesthesiology},
  author={F{\"u}gener, Andreas and Brunner, Jens O and Podtschaske, Armin},
  journal={International Journal of Production Research},
  volume={53},
  number={24},
  pages={7465--7487},
  year={2015}
}

@article{carter2001scheduling,
  title={Scheduling emergency room physicians},
  author={Carter, Michael W and Lapierre, Sophie D},
  journal={Health Care Management Science},
  volume={4},
  pages={347--360},
  year={2001}
}

@article{brunner2011long,
  title={Long term staff scheduling of physicians with different experience levels in hospitals using column generation},
  author={Brunner, Jens O and Edenharter, G{\"u}nther M},
  journal={Health Care Management Science},
  volume={14},
  pages={189--202},
  year={2011}
}

@article{ozder2020systematic,
  title={A systematic literature review for personnel scheduling problems},
  author={{\"O}zder, Emir H{\"u}seyin and {\"O}zcan, Evrencan and Eren, Tamer},
  journal={International Journal of Information Technology \& Decision Making},
  volume={19},
  number={06},
  pages={1695--1735},
  year={2020}
}

@article{van2013personnel,
  title={Personnel scheduling: A literature review},
  author={{Van den Bergh}, Jorne and Beli{\"e}n, Jeroen and De Bruecker, Philippe and Demeulemeester, Erik and De Boeck, Liesje},
  journal={European Journal of Operational Research},
  volume={226},
  number={3},
  pages={367--385},
  year={2013}
}

@article{erhard2018state,
  title={State of the art in physician scheduling},
  author={Erhard, Melanie and Schoenfelder, Jan and F{\"u}gener, Andreas and Brunner, Jens O},
  journal={European Journal of Operational Research},
  volume={265},
  number={1},
  pages={1--18},
  year={2018},
}

@Article{Bodenheimer+Smith:primary-care,
  author =       {T. S. Bodenheimer and M. D. Smith},
  title =        {Primary Care: Proposed Solutions To The Physician Shortage Without Training More Physicians},
  journal = 	 	 {Health Affairs},
  year = 	 			 {2013},
  volume =	 		 {32},
  number =	 		 {11},
  pages =	 			 {1881--1886}
}
\clearpage
\appendix
\section{Constraints of the MIP Model}\label{PR:sec:appendix_MIP_formulation}
This section provides the mathematical formulations of the constraints of the MIP presented in Section~\ref{PR:sec:MIP}. Note that we use~$\lfloor\cdot\rceil$ to denote rounding of parameters to the nearest integer.

\medskip

\noindent
\textbf{Duty Assignment:} Mandatory / non-mandatory duties must be assigned to exactly / at most one physician, respecting existing manual assignments.
\begin{align*}
\sum_{p \in \setPhysicians} \varDutyAssignment_{p, d} &= 1, \quad \forall d \in \setMandatoryDuties \tag{\ref{PR:con:mandatory_duty}} \label{PR:app:mandatory_duty} \\
\sum_{p \in \setPhysicians} \varDutyAssignment_{p, d} &\leq 1, \quad \forall d \in \setDuties\setminus\setMandatoryDuties \tag{\ref{PR:con:max_one_duty}} \label{PR:app:max_one_duty} \\
\varDutyAssignment_{p, d} &= 1, \quad \forall p \in \setPhysicians,\; d \in \setManuallyPlannedDutyAssignments(p)
\tag{\ref{PR:con:fixed_duty}} \label{PR:app:fixed_duty} \\
\varDutyAssignment_{p, d} &= 0, \quad \forall p \in \setManuallyPlannedPhysicians ,\; d \in \setDuties\setminus\setManuallyPlannedDutyAssignments(p)
\tag{\ref{PR:con:manual_duty}} \label{PR:app:manual_duty}
\end{align*}

\medskip
\noindent
\textbf{Qualifications:} Physicians can only be assigned to duties and shifts for which they are qualified.

\noindent
\begin{align*}
\varDutyAssignment_{p, d} &= 0, \quad \forall p \in \setPhysicians,\; d \in \setDuties\setminus\setQualifiedDuties(p) \tag{\ref{PR:con:qualification_duties}.1} \label{PR:app:qualification_duties}\\
\varShiftAssignment_{p, s} &= 0, \quad \forall p \in \setPhysicians,\; s \in\setShifts\setminus\setQualifiedShift(p) \tag{\ref{PR:con:qualification_duties}.2} \label{PR:app:qualification_shifts}
\end{align*}

\medskip
\noindent
\textbf{Shift Assignment:} Each shift must satisfy minimum and maximum staffing levels and adhere to qualification and manual-planning rules.

\noindent
\begin{align*}
\sum_{p \in \setShiftPhysician(s)} \varShiftAssignment_{p, s} &\geq \parShiftMinimumNumber(s), \quad \forall s \in \setShifts \tag{\ref{PR:con:shift_1}} \label{PR:app:shift_1} \\
\sum_{p \in \setShiftPhysician(s)} \varShiftAssignment_{p, s} &\leq \parShiftMaxNumber(s), \quad \forall s \in \setShifts \tag{\ref{PR:con:shift_2}}\label{PR:app:shift_2} \\
\varShiftDesiredNumber_{s} + \varShiftMaximalNumber_{s} &= 
\sum_{p \in \setShiftPhysician(s)} \varShiftAssignment_{p, s} - \parShiftMinimumNumber(s), \quad \forall s \in \setShifts \tag{\ref{PR:con:shift_3}} \label{PR:app:shift_3} \\
\varShiftDesiredNumber_{s} &\leq \parShiftDesiredNumber(s) - \parShiftMinimumNumber(s), \quad \forall s \in \setShifts
\tag{\ref{PR:con:upper_bound_shift_desired}}\label{PR:app:upper_bound_shift_desired} \\
\varShiftMaximalNumber_{s} &\leq \parShiftMaxNumber(s) - \parShiftDesiredNumber(s), \quad \forall s \in \setShifts \tag{\ref{PR:con:upper_bound_shift_max}} \label{PR:app:upper_bound_shift_max}
\end{align*}

\noindent
The following two constraints ensure that the variable~$\varShiftMaximalNumber_{s}$ can take a positive value only if $\varShiftDesiredNumber_{s}$ is at its upper bound. 
Here, the auxiliary variable~$\varShiftDesiredNumberAuxiliary_{s}$ equals~1 if $\varShiftDesiredNumber_{s}$ is \emph{not} at its upper bound.
\begin{align*}
(\parShiftDesiredNumber(s) - \parShiftMinimumNumber(s)) \cdot \varShiftDesiredNumberAuxiliary_{s} + \varShiftDesiredNumber_{s} &\geq \parShiftDesiredNumber(s) - \parShiftMinimumNumber(s), \quad \forall s \in \setShifts
\tag{\ref{PR:con:shift_4}.1}\label{PR:app:shift_4} \\
\varShiftMaximalNumber_{s} + \parShiftMaxNumber(s) \cdot \varShiftDesiredNumberAuxiliary_{s} &\leq \parShiftMaxNumber(s), \quad \forall s \in \setShifts \tag{\ref{PR:con:shift_4}.2} 
\end{align*}

\noindent
\begin{align*}
\varShiftAssignment_{p, s} &= 0, \quad \forall s \in \setShifts,\; p \in \setPhysicians\setminus\setShiftPhysician(s)
\tag{\ref{PR:con:shift_6}} \label{PR:app:shift_6} \\
\varShiftAssignment_{p, s} &= 1, \quad \forall p \in \setPhysicians,\; s \in \setManuallyPlannedShiftAssignments(p)\tag{\ref{PR:con:fixed_shift}} \label{PR:app:fixed_shift} \\
\varShiftAssignment_{p, s} &= 0, \quad \forall p \in \setManuallyPlannedPhysicians,\; s \in \setShifts\setminus\setManuallyPlannedShiftAssignments(p) \tag{\ref{PR:con:man_planned_physicians_exclude}} \label{PR:app:man_planned_physicians_exclude}
\end{align*}

\medskip
\noindent
\textbf{Rest Times:} Assignments must / should respect mandatory / desired rest times between consecutive duties or shifts.
\medskip

\noindent
Mandatory rest time between duties:
\begin{align*}
\varDutyAssignment_{p, d_1} + \varDutyAssignment_{p, d_2} &\leq 1, \quad\forall p\in\setPhysicians, (d_1,d_2)\in \setConflicts: d_1,d_2\in\setDuties
\tag{\ref{PR:con:hard_rest_duties}.1} \label{PR:app:hard_rest_duties}
\end{align*}

\noindent
Mandatory rest time between a duty and a shift:
\begin{align*}
\varDutyAssignment_{p, d_1} + \varShiftAssignment_{p, d_2} &\leq 1, \quad\forall p\in\setPhysicians, (d_1,d_2)\in \setConflicts: d_1\in\setDuties, d_2\in\setShifts
\tag{\ref{PR:con:hard_rest_duties}.2} \label{PR:app:hard_rest_duty_shift}
\end{align*}

\noindent
Mandatory rest time between a shift and a duty:
\begin{align*}
\varShiftAssignment_{p, d_1} + \varDutyAssignment_{p, d_2} &\leq 1, \quad\forall p\in\setPhysicians, (d_1,d_2)\in \setConflicts: d_1\in\setShifts, d_2\in\setDuties
\tag{\ref{PR:con:hard_rest_duties}.3} \label{PR:app:hard_rest_shift_duty}
\end{align*}

\noindent
Mandatory rest time between shifts:
\begin{align*}
\varShiftAssignment_{p, d_1} + \varShiftAssignment_{p, d_2} &\leq 1, \quad\forall p\in\setPhysicians, (d_1,d_2)\in \setConflicts: d_1, d_2\in\setShifts
\tag{\ref{PR:con:hard_rest_duties}.4} \label{PR:app:hard_rest_shifts}
\end{align*}

\noindent
Desired rest time between duties:
\begin{align*}
\varDutyAssignment_{p, d_1} + \varDutyAssignment_{p, d_2} - \varVioRestTime_{d_1, d_2}  &\leq 1, \quad\forall p\in\setPhysicians, (d_1,d_2)\in \setConflictsSoft: d_1,d_2\in\setDuties
\tag{\ref{PR:con:soft_rest_duties}.1} \label{PR:app:soft_rest_duties}
\end{align*}

\noindent
Desired rest time between a duty and a shift:
\begin{align*}
\varDutyAssignment_{p, d_1} + \varShiftAssignment_{p, d_2} - \varVioRestTime_{d_1, d_2} &\leq 1, \quad\forall p\in\setPhysicians, (d_1,d_2)\in \setConflictsSoft: d_1\in\setDuties, d_2\in\setShifts
\tag{\ref{PR:con:soft_rest_duties}.2} \label{PR:app:soft_rest_duty_shift}
\end{align*}

\noindent
Desired rest time between a shift and a duty:
\begin{align*}
\varShiftAssignment_{p, d_1} + \varDutyAssignment_{p, d_2} - \varVioRestTime_{d_1, d_2} &\leq 1, \quad\forall p\in\setPhysicians, (d_1,d_2)\in \setConflictsSoft: d_1\in\setShifts, d_2\in\setDuties
\tag{\ref{PR:con:soft_rest_duties}.3} \label{PR:app:soft_rest_shift_duty}
\end{align*}

\noindent
Desired rest time between shifts:
\begin{align*}
\varShiftAssignment_{p, d_1} + \varShiftAssignment_{p, d_2} - \varVioRestTime_{d_1, d_2} &\leq 1, \quad\forall p\in\setPhysicians, (d_1,d_2)\in \setConflictsSoft: d_1, d_2\in\setShifts
\tag{\ref{PR:con:soft_rest_duties}.4} \label{PR:app:soft_rest_shifts}
\end{align*}

\medskip
\noindent
\textbf{Absences:} Physicians cannot be assigned to duties or shifts during their absences, or to specific duties immediately before or after them.
\begin{align*}
\sum_{d \in \setDuties(t)} \varDutyAssignment_{p, d} &= 0, \quad \forall p \in \setPhysicians,\; t \in \setAbsences(p) \tag{\ref{PR:con:absence_duty}.1} \label{PR:app:absence_duty} \\
\sum_{s \in \setShifts(t)} \varShiftAssignment_{p, s} &= 0, \quad \forall p \in \setPhysicians,\; t \in \setAbsences(p) \tag{\ref{PR:con:absence_duty}.2} \label{PR:app:absence_shift} \\
\varDutyAssignment_{p, d} &= 0, \quad \forall p \in \setPhysicians,\; d \in \setDutyAbsences(p) \tag{\ref{PR:con:absence_duty_specific}.1} \label{PR:app:absence_duty_specific} \\
\varShiftAssignment_{p, s} &= 0, \quad \forall p \in \setPhysicians,\; s \in \setShiftAbsences(p) \tag{\ref{PR:con:absence_duty_specific}.2}\label{PR:app:absence_shift_specific} \\
\sum_{d \in \setDutiesBeforeAbsence\cap\setDuties(t-1)} \varDutyAssignment_{p, d} &= 0, \quad \forall p \in \setPhysicians,\; t \in \setAbsences(p) \tag{\ref{PR:con:duties_before_absence}} \label{PR:app:duties_before_absence} \\
\sum_{d \in \setDutiesAfterAbsence\cap\setDuties(t+1)} \varDutyAssignment_{p, d} &= 0, \quad \forall p \in \setPhysicians,\; t \in \setAbsences(p) \tag{\ref{PR:con:duties_after_absence}} \label{PR:app:duties_after_absence} \\
\end{align*}

\medskip
\noindent
\textbf{Duty and Shift Blocks:} All duties or shifts within a block must be assigned consistently, and required free days after blocks must be respected.
\begin{align*}
\varDutyAssignment_{p, d} &= \varDutyBlockAssignment_{p, b}, \quad \forall p \in \setPhysicians,\; b \in \setDutyBlocks: d \in b \tag{\ref{PR:con:duty_block_assignment}}\label{PR:app:duty_block_assignment} \\
\varShiftAssignment_{p, s} &= \varShiftBlockAssignment_{p, b}, \quad \forall p \in \setPhysicians,\; b \in \setShiftBlocks: s \in b \tag{\ref{PR:con:shift_block_assignment}} \label{PR:app:shift_block_assignment}
\end{align*}
\noindent
Free days after duty blocks:
\begin{align*}
\varDutyAssignment_{p, d} + \varDutyBlockAssignment_{b} \leq 1 \quad &\forall p \in \setPhysicians,\; b \in \setDutyBlocks,\; \delta \in \{1, \dots, \parRequiredFreeDaysAfterDutyBlock(b)\}, \\ &d\in\setDuties(\min\{\parEnd(b)+\delta, \parEndOfDays \})\tag{\ref{PR:con:free_days_after_duty_block_assignment}.1} \label{PR:app:free_days_after_duty_block_assignment}
\end{align*}

\begin{align*}
\varShiftAssignment_{p, s} + \varDutyBlockAssignment_{b} \leq 1 \quad &\forall p \in \setPhysicians,\; b \in \setDutyBlocks,\; \delta \in \{1, \dots, \parRequiredFreeDaysAfterDutyBlock(b)\}, \\ &s \in\setShifts(\min\{\parEnd(b)+\delta, \parEndOfDays \})\tag{\ref{PR:con:free_days_after_duty_block_assignment}.2}
\end{align*}

\noindent
Free days after shift blocks:
\begin{align*}
\varDutyAssignment_{p, d} + \varShiftBlockAssignment_{b} \leq 1 \quad &\forall p \in \setPhysicians,\; b \in \setShiftBlocks,\; \delta \in \{1, \dots, \parRequiredFreeDaysAfterDutyBlock(b)\}, \\&d\in\setDuties(\min\{\parEnd(b)+\delta, \parEndOfDays \})\tag{\ref{PR:con:free_days_after_duty_block_assignment}.3} \label{PR:app:free_days_after_shift_block_assignment}
\end{align*}

\begin{align*}
\varShiftAssignment_{p, s} + \varShiftBlockAssignment_{b} \leq 1 \quad &\forall p \in \setPhysicians,\; b \in \setShiftBlocks,\; \delta \in \{1, \dots, \parRequiredFreeDaysAfterDutyBlock(b)\}, \\ &s\in\setShifts(\min\{\parEnd(b)+\delta, \parEndOfDays \})\tag{\ref{PR:con:free_days_after_duty_block_assignment}.4}
\end{align*}

\noindent
No additional duties or shifts during blocks for which this is not allowed:
\begin{align*}
\varDutyAssignment_{p, d} + \varDutyBlockAssignment_{b} \leq 1 \quad &\forall p \in \setPhysicians,\; b \in \setDutyBlocks\cap\setBlocksNoAddDuty,\; \delta \in \{0, \dots, \parEnd(b) - \parStart(b)\}, \\ &d\in\setDuties(\min\{\parStart(b)+\delta, \parEndOfDays \})\setminus \{b\}\tag{\ref{PR:con:no_other_duties_in_duty_block}.1} \label{PR:app:no_other_duties_in_duty_block}
\end{align*}

\begin{align*}
\varShiftAssignment_{p, s} + \varDutyBlockAssignment_{b} \leq 1 \quad & \forall p \in \setPhysicians,\; b \in \setDutyBlocks\cap\setBlocksNoAddShift,\; \delta \in \{0, \dots, \parEnd(b) - \parStart(b)\}, \\&s \in\setShifts(\min\{\parStart(b)+\delta, \parEndOfDays \})\tag{\ref{PR:con:no_other_duties_in_duty_block}.2}
\end{align*}

\noindent

\begin{align*}
\varDutyAssignment_{p, d} + \varShiftBlockAssignment_{b} \leq 1 \quad &\forall p \in \setPhysicians,\; b \in \setShiftBlocks\cap\setBlocksNoAddDuty,\; \delta \in \{0, \dots, \parEnd(b) - \parStart(b)\}, \\ &d\in\setDuties(\min\{\parStart(b)+\delta, \parEndOfDays \})\tag{\ref{PR:con:no_other_duties_in_duty_block}.3} \label{PR:app:no_other_shifts_in_shift_block}
\end{align*}

\begin{align*}
\varShiftAssignment_{p, s} + \varShiftBlockAssignment_{b} \leq 1 \quad &\forall p \in \setPhysicians,\; b \in \setShiftBlocks\cap\setBlocksNoAddShift,\; \delta \in \{0, \dots, \parEnd(b) - \parStart(b)\}, \\ &s\in\setShifts(\min\{\parStart(b)+\delta, \parEndOfDays \})\setminus \{b\}\tag{\ref{PR:con:no_other_duties_in_duty_block}.4}
\end{align*}

\medskip
\noindent
\textbf{Consecutive Assignments:} Assignments across consecutive duties and blocks should comply with continuity and maximum-length requirements.
\begin{align*}
\varShiftBlockConsecPeriods_{p, b} - \varShiftBlockAssignment_{p, b} &\leq 0, \quad &\forall p \in \setPhysicians,\; b \in \setShiftBlocks \tag{\ref{PR:con:shift_block_consec_periods}.1}\label{PR:app:shift_block_consec_periods}\\
\varShiftBlockConsecPeriods_{p, b} - \varShiftBlockAssignment_{p, \parPrev(b)} &\leq 0, \quad &\forall p \in \setPhysicians,\; b \in \setShiftBlocks: \parPrev(b)\neq \texttt{None}\tag{\ref{PR:con:shift_block_consec_periods}.2}\label{PR:app:shift_block_consec_periods}\\
\varShiftBlockConsecPeriods_{p, b} &= 0, \quad &\forall p \in \setPhysicians,\; b \in \setShiftBlocks: \setPhysicianLastSB(b) \neq \emptyset, \\ & &p \in \setPhysicians\setminus\setPhysicianLastSB(b)\tag{\ref{PR:con:shift_block_consec_periods}.3}\label{PR:app:shift_block_consec_periods}
\end{align*}

\noindent
\begin{align*}
\sum_{b\in\setConsShiftBlocks_{j}}\varShiftBlockAssignment_{p, b} - \varMaxConsSBViolation_j &\leq \lvert\setConsShiftBlocks_{j}\rvert - 1, \quad\forall p \in \setPhysicians,\; j \in\{1, \dots , \lvert\setConsShiftBlocks\rvert\} \tag{\ref{PR:con:max_consecutive_shift_blocks}}\label{PR:app:max_consecutive_shift_blocks}
\end{align*}
\begin{align*}
\varDutyConsecPeriods_{p, d} - \varDutyAssignment_{p, d} &\leq 0, \quad \forall p \in \setPhysicians,\; d \in \setDuties \tag{\ref{PR:con:consecutive_duties_same_physician}.1} \label{PR:app:consecutive_duties_same_physician}\\
\varDutyConsecPeriods_{p, d} - \varDutyAssignment_{p, \parPrev(d)} &\leq 0, \quad \forall p \in \setPhysicians,\; d \in \setDuties : \parPrev(d) \neq \texttt{None} \tag{\ref{PR:con:consecutive_duties_same_physician}.2} \label{PR:app:consecutive_duties_same_physician}\\
\varDutyConsecPeriods_{p, d} &= 0, \quad \forall d \in \setDuties: \parPhysicianLastD(d) \neq\texttt{None}, \\ &\quad\quad\quad\quad p \in \setPhysicians\setminus\{\parPhysicianLastD(d)\}
\tag{\ref{PR:con:consecutive_duties_same_physician}.3} \label{PR:app:consecutive_duties_same_physician}
\end{align*}

\medskip
\noindent
\textbf{Pools:} Duties in pools must respect individual and collective limits and ensure a fair distribution.

\noindent
\begin{align*}
\sum_{d \in\setDutiesOfPool(\pi)} \varDutyAssignment_{p, d} &= \parExactDutiesPoolPerMonth(\pi), \quad \forall \pi\in\setPools,\; p\in\setPhysicianPool(\pi): \parExactDutiesPoolPerMonth(\pi)\neq\texttt{None} 
\tag{\ref{PR:con:exact_number_duties_pool}}\label{PR:app:exact_number_duties_pool}
\end{align*}
%
\begin{align*}
\sum_{d \in\setDutiesOfPool(\pi)} \varDutyAssignment_{p, d} &\leq \parMaxNumberOfDutiesInPool(\pi), \quad \forall \pi\in\setPools,\; p\in\setPhysicianPool(\pi): \parMaxNumberOfDutiesInPool(\pi)\neq\texttt{None}
\tag{\ref{PR:con:max_number_duties_pool}}\label{PR:app:max_number_duties_pool}
\end{align*}
%
\begin{align*}
\sum_{d \in\setDutiesOfPool(\pi)} \varDutyAssignment_{p, d} -\varMaxNumberOfDutiesViolation_{p, \pi}\leq \parDesiredMaxNumberOfDutiesInPool(\pi),\quad \forall \pi\in\setPools,\; p\in\setPhysicianPool(\pi):\\
\parDesiredMaxNumberOfDutiesInPool(\pi)\neq\texttt{None}
\tag{\ref{PR:con:desired_max_number_duties_pool}}\label{PR:app:desired_max_number_duties_pool}
\end{align*}

\noindent
\begin{align*}
\sum_{d \in\setDutiesOfPool(\pi)} \varDutyAssignment_{p, d} &\geq \parMinNumberOfDutiesInPool(\pi), \quad \forall \pi\in\setPools,\; p\in\setPhysicianPool(\pi): \parMinNumberOfDutiesInPool(\pi)\neq\texttt{None}
\tag{\ref{PR:con:min_number_duties_pool}}\label{PR:app:min_number_duties_pool}
\end{align*}

\noindent
\begin{align*}
\sum_{d \in\setDutiesOfPool(\pi)} \varDutyAssignment_{p, d} +\varMinNumberOfDutiesViolation_{p, \pi} \geq \parDesiredMinNumberOfDutiesInPool(\pi), \forall \pi\in\setPools,\; p\in\setPhysicianPool(\pi):\\ \parDesiredMinNumberOfDutiesInPool(\pi)\neq \texttt{None}
\tag{\ref{PR:con:desired_min_number_duties_pool}}\label{PR:app:desired_min_number_duties_pool}
\end{align*}

\noindent
\begin{align*}
\sum_{d \in\setDuties(t)\cap\setDutiesOfPool(\pi)} \sum_{p\in\setPhysicianPool(\pi)}\varDutyAssignment_{p, d} &\leq \parMaxPhysiciansPoolPerDay(\pi), \quad \forall \pi\in\setPools,\; t\in\setDays: \parMaxPhysiciansPoolPerDay(\pi)\neq\texttt{None}
\tag{\ref{PR:con:max_physicians_pool_day}}\label{PR:app:max_physicians_pool_day}
\end{align*}

\noindent
\begin{align*}
\sum_{d \in\setDuties(t)\cap\setDutiesOfPool(\pi)} \sum_{p\in\setPhysicianPool(\pi)}\varDutyAssignment_{p, d} - \varMaxNumberOfPhysiciansViolation_{p, t}\leq \parDesiredMaxPhysiciansPoolPerDay(\pi), \quad \forall \pi\in\setPools,\; t\in\setDays:\\ \parDesiredMaxPhysiciansPoolPerDay(\pi)\neq\texttt{None}
\tag{\ref{PR:con:desired_max_physicians_pool_day}}\label{PR:app:desired_max_physicians_pool_day}
\end{align*}

\noindent
\begin{align*}
\sum_{d \in\setDutiesOfPool(\pi)} \varDutyAssignment_{p, d} + \varDeviationDutiesPoolNeg_{p, \pi}\geq \lfloor\parFairDistribution(\pi, p)\rfloor, \quad \forall \pi\in\setPools,\; p\in\setPhysicianPool(\pi):\\
\parFairDistribution(\pi, p)\neq\texttt{None}
\tag{\ref{PR:con:fair_distribution_pool_min}}\label{PR:app:fair_distribution_pool_min}
\end{align*}
\begin{align*}
\sum_{d \in\setDutiesOfPool(\pi)} \varDutyAssignment_{p, d} - \varDeviationDutiesPoolPos_{p, \pi}\leq \lceil\parFairDistribution(\pi, p)\rceil, \quad \forall \pi\in\setPools,\; p\in\setPhysicianPool(\pi): \\
\parFairDistribution(\pi, p)\neq\texttt{None}
\tag{\ref{PR:con:fair_distribution_pool_max}}\label{PR:app:fair_distribution_pool_max}
\end{align*}

\medskip
\noindent
\textbf{Weekends:} Weekend assignments must / should respect hard / soft upper and lower bounds, and preferences of physicians for one/ multiple duties per weekend should be respected.
\begin{align*}
\varWeekendAttendence_{p, w} - \varDutyAssignment_{p, d} &\geq 0, \quad \forall p\in\setPhysicians,\; w\in\setWeekend,\; d\in\setDutiesOfWeekend(w) \tag{\ref{PR:con:weekend_attendance_duties}.1} \label{PR:app:weekend_attendance_duties} \\
\varWeekendAttendence_{p, w} - \sum_{d\in\setDutiesOfWeekend(w)} \varDutyAssignment_{p, d} &\leq 0, \quad \forall p\in\setPhysicians,\; w\in\setWeekend
\tag{\ref{PR:con:weekend_attendance_duties}.2} \label{PR:app:weekend_no_attendance}
\end{align*}

\noindent
\begin{align*}
\sum_{w\in\{i,\dots, i + \parDesiredMaximumConsecutiveWeekends\}\subseteq\setWeekend} \varWeekendAttendence_{p, w} &\leq \parDesiredMaximumConsecutiveWeekends, \quad \forall p\in\setPhysicians,\; i\in\setWeekend: i\leq W - \parDesiredMaximumConsecutiveWeekends
\tag{\ref{PR:con:max_consecutive_weekends}.1} \label{PR:app:max_consecutive_weekends}
\end{align*}
\begin{align*}
\sum_{w\in\{1,\dots,\parDesiredMaximumConsecutiveWeekends - \parPastConsecWeekends(p)+1\}\subseteq\setWeekend} \varWeekendAttendence_{p, w} &\leq \parDesiredMaximumConsecutiveWeekends - \parPastConsecWeekends(p), \quad \forall p\in\setPhysicians
\tag{\ref{PR:con:max_consecutive_weekends}.2} \label{PR:app:max_consecutive_weekends_2}
\end{align*}

\noindent
\begin{align*}
\varWeekendAttendence_{p, w} + \varWeekendViolation_{p, w} -\sum_{d\in\setDutiesOfWeekend(w)}\varDutyAssignment_{p, d}  &\geq 0, \quad \forall p\in\setPhysicianOneDutyWeekend,\; w\in\setWeekend
\tag{\ref{PR:con:weekend_preferences_single}} \label{PR:app:weekend_preferences_single}
\end{align*}
\begin{align*}
2 \cdot \varWeekendAttendence_{p, w} - \varWeekendViolation_{p, w} -\sum_{d\in\setDutiesOfWeekend(w)}\varDutyAssignment_{p, d} &\leq 0, \quad \forall p\in\setPhysicianSeveralDutiesWeekend,\; w\in\setWeekend
\tag{\ref{PR:con:weekend_preferences_multiple}} \label{PR:app:weekend_preferences_multiple}
\end{align*}

\noindent
\begin{align*}
\sum_{w\in\setWeekendInMonth(m)}\varWeekendAttendence_{p, w} &\leq \lfloor\parMaximumWeekends \cdot \MonthlyWeekendFactor(m)\rceil, \quad \forall p\in\setPhysicians,\; m\in\setMonths
\tag{\ref{PR:con:max_weekends_per_month}}\label{PR:app:max_weekends_per_month}
\end{align*}

\noindent
\begin{align*}
\sum_{w\in\setWeekendInMonth(m)}\varWeekendAttendence_{p, w} - \varMaxWeekendViolation_{p, m} \leq & \lfloor\parDesiredMaximumWeekends \cdot \MonthlyWeekendFactor(m)\rceil, \quad \\ & \forall p\in\setPhysicians,\; m\in\setMonths
\tag{\ref{PR:con:desired_max_weekends_per_month}}\label{PR:app:desired_max_weekends_per_month}
\end{align*}

\noindent
\begin{align*}
\sum_{w\in\setWeekendInMonth(m)}\varWeekendAttendence_{p, w} \leq & \lfloor\lvert\setWeekendInMonth(m)\rvert - \parMinFreeWeekends \cdot \MonthlyWeekendFactor(m)\rceil, \quad \\ & \forall p\in\setPhysicians,\; m\in\setMonths
\tag{\ref{PR:con:min_free_weekends_per_month}}\label{PR:app:min_free_weekends_per_month}
\end{align*}

\noindent
\begin{align*}
\sum_{w\in\setWeekendInMonth(m)}\varWeekendAttendence_{p, w} - \varMinFreeWeekendViolation_{p, m} \leq & \lfloor\lvert\setWeekendInMonth(m)\lvert - \parDesiredMinFreeWeekends \cdot \MonthlyWeekendFactor(m)\rceil, \quad \\ & \forall p\in\setPhysicians,\; m\in\setMonths
\tag{\ref{PR:con:desired_min_free_weekends_per_month}}\label{PR:app:desired_min_free_weekends_per_month}
\end{align*}

\medskip
\noindent
\textbf{Previous Planning Period:} Assignments must respect rest times following the previous planning period.
\begin{align*}
\varDutyAssignment_{p, d} &= 0, \quad \forall p \in \setPhysicians,\; d \in \setDutyRestLastPP(p) \tag{\ref{PR:con:previous_duty_rest}.1} \label{PR:app:previous_duty_rest} \\
\varShiftAssignment_{p, s} &= 0, \quad \forall p \in \setPhysicians,\; s \in \setShiftRestLastPP(p) \tag{\ref{PR:con:previous_duty_rest}.2}\label{PR:app:previous_shift_rest}\\
\end{align*}

\end{document}